\documentclass{siamltex}
\usepackage[english]{babel}
\usepackage{amsmath}
\usepackage{amsfonts}
\usepackage{amssymb}
\usepackage{graphicx}
 \usepackage{array,multirow,makecell}
\usepackage{tikz}
\usepackage{here}
\usepackage{algorithm}
\usepackage{algorithmic}
\usepackage{graphicx}
\usepackage{amsfonts}            
\usepackage{amsmath}
\usepackage{here}
\usetikzlibrary{matrix,calc}
\usepackage{caption}

\newtheorem{lemm}{Lemma}
\newtheorem{pf}{Proof}

\newtheorem{theo}{Theorem}

\newtheorem{defis}{Definitions}
\newtheorem{defi}[theo]{Definition}
   
\title{On  the tensor nuclear norm and the total variation regularization for  image and video completion}

\author{ A.H Bentbib\thanks{Facult\'e des Sciences et Techniques-Gueliz, Laboratoire de Math\'ematiques Appliqu\'ees et Informatique, Marrakech, Morocco} \and A. El Hachimi\footnotemark[4] 
	\and K. Jbilou\footnotemark[4] \thanks{Universit\'e du Littoral Cote d'Opale, LMPA, 50 rue F. Buisson, 62228 Calais-Cedex, France }
	\and {A. Ratnani  \thanks{laboratory MSDA, Mohammed VI Polytechnic University, Green City, Morocco}.}}
\begin{document}
	
	\maketitle

	\begin{abstract}
		In the present  paper we  propose two new algorithms of tensor completion for three-order  tensors. The proposed methods consist   in minimizing the average rank of the underlying tensor using its approximate function namely the tensor nuclear norm and then  the  recovered data will be obtained  by using the total variation regularisation technique. We will adopt the Alternating Direction Method of Multipliers (ADM), using the tensor T-product,  to solve the main optimization problems associated  to  the two algorithms. In the last section, we present some numerical experiments and comparisons with the most known image  completion methods.
	\end{abstract}
	
	
	\section{Introduction}
	Tensors become an important notion that treat the high  dimensional data and it plays an important role in a wide range of real-world applications.
	In  this work, we will be interested in  the problem of tensor completion, with the aim of estimating the missing values from an observed data, e.g.,   inpainting color images \cite{bertalmio2000image,chan2011augmented,gao2018mixture,komodakis2006image}, hyperspectral image recovery \cite{chen2018destriping,elguide,jiang2018matrix,li2012coupled}, magnetic resonance image recovery \cite{varghees2012adaptive,fu20163d}, and higher order web site link analysis \cite{kolda2005higher}.\\
	The main idea behind the completion problem is to find a low-rank data containing the main information of the original data. For matrices,  the problem of completion is mathematically formulated as follows
	\begin{eqnarray}\label{eq1}
	&\underset{A}{min}&\; rank(A)\nonumber \\
	&s.t& P_\Omega (A)=P_\Omega(M),
	\label{eq 1}
	\end{eqnarray} 
	where $A\in \mathbb{R}^{n_1\times n_2}$ is the underlying matrix, $M\in \mathbb{R}^{n_1\times n_2}$ is the observed matrix, $\Omega$ is the set of known pixels and $P_\Omega$ is the projection operator onto $\Omega$.  This optimization problem is not easy to solve because of the non-convexity of the rank function . For that reason, Fazel \cite{fazel2001rank}, and Kurucz \cite{kurucz2007methods}  proposed to solve the problem \eqref{eq1}  by using the convex surrogate of the rank and  the SVD decomposition.  In \cite{candes2008exact}, Candes and Recht  proved theoretically  that under certain conditions,  the following optimization problem 
	\begin{eqnarray}\label{eq2}
	&\underset{A}{min}&\; \left\Vert A\right\Vert_*\nonumber \\
	&s.t& P_\Omega (A)=P_\Omega(M)
	\end{eqnarray}
	recovers well the data, where $\left\Vert A\right\Vert_*$ is the nuclear norm of $A$ which will be defined later. Since tensors are the generalization of matrices,  the problem of tensor completion can  be formulated as follows
	\begin{eqnarray}\label{eq3}
	&\underset{\mathcal{A}}{min}&\; rank(\mathcal{A})\nonumber\\
	& s.t & \mathcal{P}_\Omega (\mathcal{A})=\mathcal{P}_\Omega (\mathcal{M})
	\end{eqnarray}
	where $\mathcal{A}$ is the underlying tensor of order 3, $\mathcal{M}$ is the observed tensor, $\Omega$ is the set of the known data and $\mathcal{P}_\Omega $ is the projection operator defined by\\
	\begin{equation}\label{eq4}
	\mathcal{P}_{\Omega}(\mathcal{A})_{i,j,k}=\left\{\begin{array}{rl}
	\mathcal{A}_{i,j,k} & \; , \; (i,j,k)\in \Omega \\
	0 &\; \; otherwise.
	\end{array}\right. \nonumber
	\end{equation}
	As stated in \cite{hillar2013most}, the optimization problem \eqref{eq3} is  N-P hard and then one should study the tensor version of the problem \eqref{eq2} and this is the main subject of the present work.   The new tensor-rank optimisation problem will be solved using the tensor T-product which is based on the Fast Fourier Transform (FFT).  Notice that  the notion of tensor rank is complicated as compared to the matrix rank and  many tensor rank definitions and procedures such as  Tucker-rank \cite{kolda2009tensor}, CP-rank \cite{kolda2009tensor}, TT-rank \cite{oseledets2011tensor} and the  tensor tubal rank \cite{kilmer2013third}, have been introduced the last years;  see also  \cite{ding2019low,ji2016tensor,bengua2017efficient,xu2013parallel}.\\
	
	\noindent The outline of this paper is as follows:
	In Section \ref{sec 2} we give some notations and preliminaries that will be used in the paper.   Section \ref{sec 3} is devoted to the development of  our proposed tensor completion approaches. We will show how to use the tensor nuclear norm in combination with the TV regularisation procedure to derive the  new  completion algorithms. In the two approaches, we will use the well known tensor T-product. Some numerical experiments with comparisons to the most well known  methods are   presented in Section  \ref{sec 4}, showing the effectiveness of the presented approaches.

	\begin{section}{Notations and preliminaries}
		\label{sec 2}
		In this paper we denote tensors by calligraphic letters,
		e.g., $\mathcal{A}$. Matrices are denoted by  capital letters, e.g.,
		$A$ and  vectors are denoted by  lower case letters, e.g., $a$.\\
		Let $\mathcal{A}\in \mathbb{R}^{n_1\times n_2\times n_3}$ be an $3$-order tensor.  We define its  Frobenius-norm by
		\begin{equation*}
		\left\Vert \mathcal{A}\right\Vert_F=\sqrt{\sum_{i_1=1}^{n_1}\sum_{i_2=1}^{n_2}\sum_{i_3=1}^{n_3} a_{i_1,i_2,i_3}^2}.
		\end{equation*}
		The inner product between the two tensors $\mathcal{A}$ and $\mathcal{B}$ in $\mathbb{R}^{n_1\times n_2\times n_3}$ is given by 
		\begin{equation*}
		\left<\mathcal{A},\mathcal{B}\right>=\sum_{i_1=1}^{n_1}\sum_{i_2=1}^{n_2}\sum_{i_3=1}^{n_3} a_{i_1,i_2,i_3} b_{i_1,i_2,i_3}.
		\end{equation*}
		Let $\mathcal{A}\in \mathbb{R}^{n_1\times n_2\times n_3}$, then the  $i^{th}$ frontal slice of the tensor  $\mathcal{A}$  is  the matrix $\mathcal{A}(:,:,i)$ and will be denoted by  $\mathcal{A}^{(i)}$. 
		For two matrices $A\in \mathbb{R}^{n\times m}$ and $B\in \mathbb{R}^{p\times q}$, the  Kronecker product is the $np\times mq$ matrix  given as 
		\begin{equation*}
		A\otimes B =[a_{ij}B]_{i=1:n;j=1:m}
		\in \mathbb{R}^{np\times mq}.
		\end{equation*}
		\begin{subsection}{The discrete Fourier transform}
			Let $v \in \mathbb{R}^n$, we denote its Discrete Fourier Transform (DFT) by $\hat{v}$ and it is defined by 
			\begin{equation*}
			\hat{v}=F_n v
			\label{eq 4}
			\end{equation*}
			where $F_n$ denotes the DFT matrix and it is defined by 
			\begin{equation}\label{four1}
			F_n=\begin{pmatrix}
			1 & 1 & 1 & \dots & \dots & 1 & 1 \\
			1 & \omega & \omega^2 & \omega^3 & \dots & \omega^{n-1} & \omega^{n}\\
			\vdots & \vdots & \vdots & \vdots & \vdots & \vdots & \vdots  \\ 
			1 & \omega^{n-1} & \omega^{2(n-1} & \dots & \dots & \omega^{(n-2)(n-1)} & \omega^{(n-1)^2}
			\end{pmatrix}\in \mathbb{R}^{n\times n}
			\end{equation}
			notice that $\dfrac{F_n}{\sqrt{n}}$ is a unitary matrix, i.e.,
			\begin{equation*}
			F_n^* F_n=F_n F_n^*=nI_n.
			\label{eq 6}
			\end{equation*}
			The Fast Fourier Transform (FFT) allows us to compute the matrix-vector $F_nv$ in a very economical way. It computes this product  with a cost of  $O\left(nlogn\right)$ instead of $O\left(n^2\right)$ and it is represented in Matlab by the command $fft$ and its inverse by $ifft$:  $\hat{v}={\tt fft}(v)$ and $v={\tt ifft} (\hat{v})$.\\
			The circulant matrix associated to the vector $v$ is given as 
			\begin{equation*}
			{\tt circ}(v)=\begin{pmatrix}
			v_1 & v_n & v_{n-1} & \dots & v_2 \\
			v_2 & v_1 & v_n & \dots & v_3 \\
			\vdots & \vdots & \ddots & \dots & \vdots \\
			v_n & v_{n-1} & \dots & v_2 & v_1
			\end{pmatrix}\in \mathbb{R}^{n\times n}
			\end{equation*}
			which it can be diagonalized by using the DFT and we get
			\begin{equation*}
			F_n {\tt circ} (v)F_n^{-1}={\tt diag}(\hat{v})
			\label{eq 8}
			\end{equation*}
			with ${\tt diag} (\hat{v})$ denotes the diagonal matrix, where the $i^{th}$ element of its diagonal is $\hat{v}_i$.
			\begin{lemm} \cite{rojo2004some}
				Given a real vector $v \in \mathbb{R}^n$, the associated $\hat{v}=F_n v$ satisfies
				\begin{equation}\label{eq 9}
				\hat{v}_1\in \mathbb{R} \; {\rm and} \; {\tt conj}(\hat{v}_i)=\hat{v}_{n-i+2},  \; \; i=2,...,\left[\dfrac{n+1}{2}\right].
				\end{equation}
				Conversely, for any given complex $\hat{v}\in \mathbb{C}^n$ satisfying \eqref{eq 9}, there exists a real circulant matrix ${\tt circ}(v)$ such that \eqref{eq 8} holds.
				\label{lm 1}
			\end{lemm}
			\medskip
			
			\noindent 	Let $\mathcal{A}\in \mathbb{R}^{n_1\times n_2\times n_3}$ be a 3-order tensor, we denote its DFT along each tubes $\hat{\mathcal{A}}\in \mathbb{C}^{n_1\times n_2\times n_3},$.  This operation  can be done  in Matlab  by using the following command
			\begin{equation*}
			\hat{\mathcal{A}}={\tt fft}\left(\mathcal{A},[\,],3\right).
			\end{equation*}
			Conversely, we can obtain $\mathcal{A}$ from $\hat{\mathcal{A}}$ using the Matlab command
			\begin{equation*}
			\mathcal{A}={\tt  ifft}\left(\hat{\mathcal{A}},[\,],3\right).
			\end{equation*}
			Thanks to Lemma \ref{lm 1}, we have
			\begin{equation*}
			\hat{\mathcal{A}}^{(1)}\in \mathbb{R}^{n_1\times n_2}\;\;  and \;\; conj\left(\hat{\mathcal{A}}^{(i)}\right)=\hat{\mathcal{A}}^{(n_3 -i +2)} \;\; for \;\; i=2,...,\left[\dfrac{n_3+1}{2}\right].
			\end{equation*}
			We have also
			\begin{equation}\label{eq 13}
			\left\Vert \mathcal{A}\right\Vert_F=\dfrac{1}{\sqrt{n_3}} \left\Vert \hat{\mathcal{A}}\right\Vert_F\;\; and \;\; \left<\mathcal{A},\mathcal{B}\right>=\dfrac{1}{n_3} \left<\hat{\mathcal{A}},\hat{\mathcal{B}}\right>.
			\end{equation}
			We define the block diagonal matrix associated to the tensor  $\mathcal{A}\in \mathbb{R}^{n_1\times n_2\times n_3}$ as follows
			\begin{equation}\label{bdiag}
			{\tt 	bdiag}\left(\mathcal{A}\right)=\begin{pmatrix}
			\mathcal{A}^{(1)} & & & & \\
			&\mathcal{A}^{(2)} & & &  \\
			& & \mathcal{A}^{(3)} & &   \\
			& &  & \ddots &   \\
			& & & &  \mathcal{A}^{(n_3)}  \\
			\end{pmatrix}
			\end{equation}
			also, we define its block circulant matrix by 
			\begin{equation}\label{bcirc}
			{\tt	bcirc }\left(\mathcal{A}\right)=\begin{pmatrix}
			\mathcal{A}^{(1)} &\mathcal{A}^{(n_3)} &\dots &\dots &\mathcal{A}^{(2)} \\
			\mathcal{A}^{(2)}&\mathcal{A}^{(1)} & \dots& \dots & \mathcal{A}^{(3)} \\
			\vdots&  \ddots &\ddots & \ddots & \vdots  \\
			\vdots&  \ddots &\ddots & \ddots & \vdots  \\
			\mathcal{A}^{(n_3)} & \mathcal{A}^{(n_3 -1)} & \dots & \dots & \mathcal{A}^{(1)}.
			\end{pmatrix}
			\end{equation}
			As  $	{\tt bcirc}(\mathcal{A})$ is a block circulant matrix, it can be block diagonalized  using the DFT \cite{hao2013facial}.  Then we  get 
			\begin{equation}\label{eq 16}
			\left(F_{n_3}\otimes I_{n_1}\right) 	{\tt bcirc}(\mathcal{A}) \left(F_{n_3}^{*}\otimes  I_{n_2} \right)=	{\tt bdiag}(\hat{\mathcal{A}}).
			\end{equation}
		\end{subsection}
		\begin{subsection}{The tensor T-product}
			Let $\mathcal{A}\in \mathbb{R}^{n_1\times n_2 \times n_3}$ a third-order tensor, we define the following operators
			\begin{equation*}
			{\tt unfold}(\mathcal{A})=\left [	\mathcal{A}^{(1)},
			\mathcal{A}^{(2)}, \ldots,
			\mathcal{A}^{(n_3)} \right]^T,\;\; 	{\tt fold(unfold }(\mathcal{A}))=\mathcal{A}.
			\end{equation*}
			\begin{defi}{T-product \cite{kilmer2011factorization}}\\
				Let $\mathcal{A}\in \mathbb{R}^{n_1\times n\times n_3}$ and $\mathcal{B}\in \mathbb{R}^{n\times n_2\times n_3 }$, we define the t-product between $\mathcal{A}$ and $\mathcal{B}$ by
				\begin{equation}\label{eq 18}
				\mathcal{A}*\mathcal{B}=fold\left(	{\tt bcirc}\left(\mathcal{A}\right)	{\tt unfold}\left(B\right)\right).
				\end{equation}
			\end{defi}
			We notice that from \eqref{eq 16},  we can compute the T-product between two tensors $\mathcal{A}$ and $\mathcal{B}$ of appropriate sizes using the following property
			\begin{equation*}
			\mathcal{D}=\mathcal{A}*\mathcal{B}\Longleftrightarrow 	{\tt bdiag}(\hat{\mathcal{D}})=	{\tt  bdiag}(\hat{\mathcal{A}})	{\tt bdiag}(\hat{\mathcal{B}}).
			\end{equation*}
			The following algorithm summarises the different steps for the T-product
			\begin{algorithm}[H]
				\caption{The T-product via the  FFT.}
				\label{alg:1}
				\begin{algorithmic}[1]
					\STATE \textbf{Inputs:} $\mathcal{A}\in \mathbb{R}^{n_1\times n\times n_3}$ and $\mathcal{B}\in \mathbb{R}^{n\times n_2\times n_3}$ 
					\STATE \textbf{Output:} $\mathcal{C}=\mathcal{A}*\mathcal{B}\in \mathbb{R}^{n_1\times n_2\times n_3}$ 
					\STATE Compute $\hat{\mathcal{A}}={\tt fft}\left(\mathcal{A},[\,],3\right)$ and $\hat{\mathcal{B}}={\tt fft}\left(\mathcal{B},[\, ],3\right)$
					\FOR {$i=1$ to $\left[\dfrac{n_3+1}{2}\right]$}
					\STATE  $\hat{\mathcal{C}}(:,:,i)=\hat{\mathcal{A}}(:,:,i)\hat{\mathcal{B}}(:,:,i)$
					\ENDFOR
					\FOR {$i=\left[\dfrac{n_3+1}{2}\right]+1$ to $n_3$}
					\STATE  $\hat{\mathcal{C}}(:,:,i)=conj\left(\hat{\mathcal{C}}(:,:,n_3+2-i)\right)$
					\ENDFOR
					\STATE $\hat{\mathcal{C}}(:,:,i)=\hat{\mathcal{A}}(:,:,i)\hat{\mathcal{B}}(:,:,i)$
					\STATE $\mathcal{C}={\tt ifft}(\hat{\mathcal{C}},[\,],3)$
				\end{algorithmic}
			\end{algorithm}
			\noindent 	Using Algorithm \ref{alg:1},  the cost of computing the T-product of  $\mathcal{A}\in \mathbb{R}^{n_1\times n\times n_3}$ and $\mathcal{B}\in \mathbb{R}^{n\times n_2\times n_3}$ is $O\left(\dfrac{nn_1n_2n_3}{2}\right) $ instead of $O\left(n^2 n_1 n_2 n_3 \right) $ if we use directly  the relation \eqref{eq 18}.
		\end{subsection}
		\begin{subsection}{The tensor SVD}
			In the sequel, we need the following definitions.
			\begin{defis}\cite{kilmer2011factorization}
				\begin{itemize}
					\item \textbf{Conjugate transpose:} The conjugate transpose of a tensor $\mathcal{A}\in \mathbb{C}^{n_1\times n_2 \times n_3}$ is the tensor $\mathcal{A}^*\in\mathcal{C}^{n_2\times n_1\times n_3}$ obtained by conjugate transposing each of its frontal slices and then reversing the order of transposed frontal slices $2$ through $n$.
					\item \textbf{Identity tensor:} The identity tensor $\mathcal{I}\in \mathcal{R}^{l\times l \times n}$ is the tensor with its first frontal slice being the $l\times l$ identity matrix, and other frontal slices being all zeros.
					\item \textbf{F-diagonal tensor:} A tensor is called f-diagonal if each of its frontal slices is a diagonal matrix.
					\item \textbf{Orthogonal tensor:} A tensor $\mathcal{Q}\in \mathcal{R}^{l\times l\times n}$ is orthogonal if it satisfies 
					$$\mathcal{Q}*\mathcal{Q}^*=\mathcal{Q}^**\mathcal{Q}=\mathcal{I}.$$
				\end{itemize}
			\end{defis}
			The Singular Value Decomposition (SVD) for matrices, was generalized to the tensor case using the T-product \cite{kilmer2011factorization} as is stated in the following theorem. 
			\begin{theo}
				Let $\mathcal{A}\in \mathbb{R}^{n_1\times n_2\times n_3}$ real valued tensor, then $\mathcal{A}$ can be factored as
				$$\mathcal{A}=\mathcal{U}*\mathcal{S}*\mathcal{V}^T$$
				with $\mathcal{U}\in \mathbb{R}^{n_1\times n_1\times n_3}$ and $\mathcal{V}\in \mathcal{R}^{n_2\times n_2 \times n_3}$ are orthogonal tensors, and $\mathcal{S}\in \mathbb{R}^{n_1\times n_2\times n_3}$ is an $f$-diagonal tensor.
			\end{theo}
			
			\medskip
			\noindent The process called T-SVD of decomposing a 3-order tensor via the tensor T-product  is summarized in the following algorithm 
			\begin{algorithm}[H]
				\caption{T-SVD}
				\label{alg:3}
				\begin{algorithmic}[1]
					\STATE \textbf{Impute} $\mathcal{A}\in \mathbb{R}^{n_1\times n_2\times n_3}$
					\STATE \textbf{Output:} t-SVD components $\mathcal{U}$, $\mathcal{S}$ and $\mathcal{V}$
					\STATE  $\hat{\mathcal{A}}={\tt fft} (\mathcal{A},[\,],3)$
					\FOR {$i=1$ to $\left[\dfrac{n_3+1}{2}\right]$}
					\STATE  $[\hat{\mathcal{U}}(:,:,i),\hat{\mathcal{S}}(:,:,i),\hat{\mathcal{V}}(:,:,i)]=SVD(\hat{\mathcal{A}}(:,:,i))$
					\ENDFOR
					\FOR {$i=\left[\dfrac{n_3+1}{2}\right]+1$ to $n_3$}
					\STATE $\hat{\mathcal{U}}(:,:,i)={\tt conj}(\hat{\mathcal{U}}(:,:,n_3+2-i))$
					\STATE $\hat{\mathcal{S}}(:,:,i)={\tt conj}(\hat{\mathcal{S}}(:,:,n_3+2-i))$
					\STATE $\hat{\mathcal{V}}(:,:,i)={\tt conj}(\hat{\mathcal{V}}(:,:,n_3+2-i))$
					\ENDFOR
					\STATE $\mathcal{U}={\tt ifft}(\hat{\mathcal{U}},[\,],3)$, $\mathcal{S}={\tt  ifft}(\hat{\mathcal{S}},[\, ],3)$ and $\mathcal{V}={\tt ifft}(\hat{\mathcal{V}},[\,],3)$
				\end{algorithmic}
			\end{algorithm}
			
			\medskip
			\noindent Next, we recall the definitions of the tensor  tubal rank \cite{lu2019tensor} and the tensor average rank \cite{lu2019tensor} that will be used in this paper. 
			\begin{defi}
				For $\mathcal{A}\in \mathbb{R}^{n_1\times n_2\times n_3}$, the tensor tubal rank, denoted as ${\tt rank_t}\left(\mathcal{A}\right)$, is the number of nonzero singular tubes of $\mathcal{S}$, where $\mathcal{S}$ is from the t-SVD of $\mathcal{A}=\mathcal{U}*\mathcal{S}*\mathcal{V}^*$.
				We can write
				\begin{equation}
				{\tt rank_t}\left(\mathcal{A}\right)={\tt card}\left(\left\{i / \; \mathcal{S}\left(i,i,:\right)\neq 0\right\}\right).
				\end{equation}
			\end{defi}
			\begin{defi}
				For $\mathcal{A}\in \mathbb{R}^{n_1\times n_2\times n_3}$, the tensor average rank, denoted as ${\tt rank_a}(\mathcal{A})$, is defined as 
				\begin{equation}
				{\tt rank_a}(\mathcal{A})=\dfrac{1}{n_3}{\tt rank(bcirc}(\mathcal{A})).
				\end{equation}
			\end{defi}
		\end{subsection}
		\begin{subsection}{The tensor nuclear norm.}  We first recall the matrix nuclear norm.  
			Let $A\in \mathbb{R}^{n\times m}$ be  a matrix, then the nuclear norm denoted by $\Vert A \Vert_*$,  is defined  as the dual norm of the matrix spectral norm, i.e.
			\begin{equation}\label{eq 22}
			\Vert A \Vert_*=\underset{\left\Vert B \right\Vert \leq 1}{arg\, min}\; \left\vert \left<A,B\right>\right\vert      		
			\end{equation}
			where $\left\Vert B \right\Vert$ denotes the matrix spectral norm.
			We notice that we also have 
			\begin{equation}
			\Vert A \Vert_*=\sum_{i=1}^{r} \sigma_i
			\end{equation}
			where $\left\{\sigma_i\right\}_{i=1}^r$ are the singular values of $A$ and $r$ is the tank of $A$.
			\begin{defi}
				Let $\mathcal{A}\in \mathbb{R}^{n_1\times n_2\times n_3}$ be a 3-order tensor. Then the  tensor spectral norm of $\mathcal{A}$ is defined as
				\begin{equation}
				\left\Vert \mathcal{A}\right\Vert = \left\Vert bcirc(\mathcal{A})\right\Vert
				\end{equation}
			\end{defi}
			Also, from \eqref{bdiag} and \eqref{bcirc}, we get
			\begin{equation}
			\left\Vert \mathcal{A}\right\Vert = \left\Vert bdiag(\hat{\mathcal{A}})\right\Vert.
			\end{equation}
			The tensor nuclear norm is an extension of the matrix nuclear norm to tensors and is defined in the following definition
			\begin{defi}
				Let $\mathcal{A}\in \mathbb{R}^{n_1\times n_2\times n_3}$ be a three mode-tensor.  The tensor nuclear norm of $\mathcal{A}$  is defined as follows 
				\begin{equation}
				\left\Vert \mathcal{A}\right\Vert_*=\underset{\left\Vert \mathcal{B}\right\Vert\leq 1}{arg\, min}\; \left\vert \left<\mathcal{A},\mathcal{B}\right>\right\vert.
				\label{eq 26}
				\end{equation}
			\end{defi}
			Using \eqref{eq 13} and \eqref{eq 26}, we get  the following relations 
			\begin{equation}
			\left\Vert \mathcal{A}\right\Vert_*=\dfrac{1}{n_3}\left\Vert bcirc(\mathcal{A})\right\Vert_*=\dfrac{1}{n_3}\left\Vert bdiag(\hat{\mathcal{A}})\right\Vert_*, 
			\end{equation}
			and 
			\begin{equation}
			\left\Vert \mathcal{A}\right\Vert_*=\sum_{i=1}^{r}\mathcal{S}(i,i,1)
			\end{equation}
			where $r$ is the tubal rank of $\mathcal{A}$, and $\mathcal{S}$ is given from the T-SVD of $\mathcal{A}$.
			\begin{theo}\label{theo 12}\cite{lu2019tensor} \\
				On the set $\left\{\mathcal{A}\in \mathbb{R}^{n_1\times n_2 \times n_3}\; / \; \left\Vert \mathcal{A}\right\Vert\leq 1 \right\}$ the convex envelope of the average rank of $\mathcal{A}$ is its tensor nuclear norm $\left\Vert \mathcal{A}\right\Vert_*$.
			\end{theo}
		\end{subsection}
		\begin{subsection}{Tensor singular value thresholding } 
			In this subsection, we first recall the Tensor singular value thresholding  \cite{lu2019tensor} and 
			give an algorithm summarizing the whole process that will be used later.
			
			\begin{defi}
				Let $\mathcal{A}\in \mathbb{R}^{n_1\times n_2\times n_3}$ be a tensor and consider its  tensor SVD decomposition as 
				$$\mathcal{A}=\mathcal{U}*\mathcal{S} *\mathcal{V}^T.$$
				The tensor singular value thresholding of $\mathcal{A}$ with a given parameter  $\tau$ is defined by
				\begin{equation}
				\mathcal{D}_\tau \left(\mathcal{A}\right)=\mathcal{U}*\mathcal{S}_{\tau} * \mathcal{V}^T
				\label{eq 24}
				\end{equation}
				where $\mathcal{S}_{\tau} =ifft\left(max\left(\hat{\mathcal{S}}-\tau,0 \right),\left[\, \right],3\right).$
			\end{defi}
			\begin{theo}\cite{lu2019tensor} \label{theo 8}
				For any $\tau >0$ and $\mathcal{A}\in \mathbb{R}^{n_1\times n_2\times n_3}$, the tensor singular value thresholding \eqref{eq 24} is connected to the nuclear norm via  the following relation 
				\begin{equation}
				\mathcal{D}_\tau \left(\mathcal{A}\right)=\underset{\mathcal{X} \in \mathbb{R}^{n_1\times n_2\times n_3}}{arg\, min}\; \tau \left\Vert \mathcal{X} \right\Vert_* + \left\Vert \mathcal{A}-\mathcal{X}\right\Vert_F^2.
				\end{equation}
			\end{theo}
			
			\medskip
			\noindent The tensor singular value thresholding  process is summarized in the following algorithm
			\begin{algorithm}[H]
				\caption{Tensor singular value thresholding algorithm}
				\label{alg:3}
				\begin{algorithmic}[1]
					\STATE \textbf{Impute} $\mathcal{A}\in \mathbb{R}^{n_1\times n_2 \times n_3}$
					\STATE \textbf{Output:} $\mathcal{D}_\tau \left(\mathcal{A}\right)$
					\STATE  $\hat{\mathcal{A}}={\tt fft}(\mathcal{A},[\, ],3)$
					\FOR {$i=1$ to $\left[\dfrac{n_3+1}{2}\right]$}
					\STATE  $[\hat{\mathcal{U}}(:,:,i),\hat{\mathcal{S}}(:,:,i),\hat{\mathcal{V}}(:,:,i)]=SVD(\hat{\mathcal{A}}(:,:,i))$
					\STATE $\hat{\mathcal{S}}^{(i)}_\tau = \left(\hat{\mathcal{S}}^{(i)}-\tau\right)_+$
					\STATE $\mathcal{D}_\tau \left(\hat{\mathcal{A}}\right)^{(i)}=\hat{\mathcal{U}}(:,:,i) \hat{\mathcal{S}}(:,:,i)_\tau \hat{\mathcal{V}}(:,:,i)$
					\ENDFOR
					\FOR {$i=\left[\dfrac{n_3+1}{2}\right]+1$ to $n_3$}
					\STATE $\mathcal{D}_\tau \left(\hat{\mathcal{A}}\right)^{(i)}=\mathcal{D}_\tau \left(\hat{\mathcal{A}}\right)^{(n_3+2-i)}$
					\ENDFOR
					\STATE $\mathcal{D}_\tau \left(\mathcal{A}\right)={\tt  ifft}(\mathcal{D}_\tau \left(\hat{\mathcal{A}}\right),[\,],3)$, 
				\end{algorithmic}
			\end{algorithm}
		\end{subsection}
	\end{section}

	\begin{section}{The proposed methods}
		\label{sec 3}
		Our proposed approaches  are based on the minimization of the average rank of a three-order tensor $\mathcal{A}\in\mathbb{R}^{n_1\times n_2 \times n_3}$. The problem  can be  formulated as follows 
		\begin{eqnarray}\label{eq 3.1}
		&\underset{\mathcal{A}}{min}&\; {\tt rank_a}\left(\mathcal{A}\right) \nonumber \\
		&s.t&  \mathcal{P}_\Omega(\mathcal{A})=\mathcal{P}_\Omega(\mathcal{M}).
		\end{eqnarray}
		where $\mathcal{A}\in \mathbb{R}^{n_1\times n_2\times n_3}$ is the underlying tensor, $\mathcal{M}\in \mathbb{R}^{n_1\times n_2\times n_3}$ is the observed tensor and  $\Omega$ is the set of the known pixels. 
		As stated in  \cite{zhang2014novel} and thanks to Theorem \ref{theo 12}, we can replace  the problem  \eqref{eq 3.1} by the following one 
		\begin{eqnarray}\label{eq 3.2}
		&\underset{\mathcal{A}}{min}&\; \left\Vert \mathcal{A}\right\Vert_* \nonumber \\
		&s.t& \mathcal{P}_\Omega(\mathcal{A})=\mathcal{P}_\Omega(\mathcal{M})
		\end{eqnarray}
		where  $\left\Vert .\right\Vert_*$ is the tensor nuclear norm defined above. 
		It is known that in  real world,   the problem \eqref{eq 3.2} can be very ill-conditioned and one needs regularisation techniques such as the   $(\textbf{TV})$-regularisation which  we will consider in the present paper. As in the matrix case, other regularization procedures are also possible.
		
		\begin{subsection}{Tensor completion using the tensor nuclear norm and the first order total variation}
			For our first approach, we consider the TV-regularized problem 
			\begin{eqnarray}
			&\underset{\mathcal{A}}{min}& \; \left\Vert \mathcal{A}\right\Vert_* + \lambda \textbf{TV}_1\left(\mathcal{A}\right) \nonumber \\
			& s.t & \mathcal{P}_\Omega (\mathcal{A})=\mathcal{P}_\Omega (\mathcal{M})
			\label{eq 3.3}
			\end{eqnarray}
			where $\lambda$ is a regularization parameter, and 
			\begin{equation}
			\textbf{TV}_1(\mathcal{A})=\left[TV_1(\mathcal{A}^{(1)})|\,TV_1(\mathcal{A}^{(2)}|\,...|\,TV_1(\mathcal{A}^{(n_3)})\right]\in \mathbb{R}^{n_1\times n_2 \times n_3}, \nonumber
			\end{equation}
			with 
			\begin{equation}
			TV_1(\mathcal{A}^{(n)})=\sum_{i=1}^{n_1}\sum_{j=1}^{n_2}\sqrt{\left(D_1^1\mathcal{A}^{(n)}\right)_{i,j}^2+\left(D_1^2\mathcal{A}^{(n)}\right)_{i,j}^2}, \; n\in \{1,2,...,n_3\} \nonumber
			\end{equation}
			and  $D_1$ and $D_2$ are  the derivative operators in the first and the second direction, respectively, with 
			$$D_1^1 \mathcal{A}^{(n)}=\mathcal{A}^{(n)} C_{n_2}^1, \;\;  \;\; D_1^2 \mathcal{A}^{(n)}=C_{n_1}^2 \mathcal{A}^{(n)} $$
			where  $C_1$ and $C_2$ are  the matrices defined as 
			$$C_m^1=\begin{pmatrix}
			-1 & 0 & 0 & \dots & 0 & 1 \\
			1 & -1 & 0 & \dots & \dots & 0 \\
			0 & 1 & -1 & 0 & \dots & \vdots \\
			\vdots & \ddots & \ddots & \ddots & \ddots & \vdots \\
			\vdots & \ddots & \ddots & \ddots & \ddots & \vdots \\
			0 & 0 & 0 & \dots & 1 & -1 
			\end{pmatrix}\in\mathbb{R}^{m\times m},\;  C_p^2=\begin{pmatrix}
			-1 & 1 & 0 & \dots & 0 & 0 \\
			0 & -1 & 1 & \dots & \dots & 0 \\
			0 & 0 & -1 & 1 & \dots & \vdots \\
			\vdots & \ddots & \ddots & \ddots & \ddots & \vdots \\
			\vdots & \ddots & \ddots & \ddots & \ddots & 1 \\
			1 & 0 & 0 & \dots & 0 & -1 
			\end{pmatrix}\in \mathbb{R}^{p\times p}.$$
			To solve the constrained optimization problem  \eqref{eq 3.3}, we have to go through the following intermediate optimization problem
			\begin{eqnarray}	\label{eq 3.4}
			&\underset{\mathcal{A},\, \mathcal{Z},\, \mathcal{Y}}{min}& \;\left [  \left\Vert\mathcal{Z}\right\Vert_*+\lambda\sum_{n=1}^{n_3}\sum_{j=1}^{n_2}\sum_{i=1}^{n_1}\left\Vert \mathcal{Y}_{i,j}^{(n)}\right\Vert_2 \right]  \\
			& s.t & \mathcal{P}_\Omega (\mathcal{A})=\mathcal{P}_\Omega(\mathcal{M}),\; \mathcal{Z}=\mathcal{A},\; \mathcal{Y}_1=\mathcal{D}_1^1\mathcal{A} \; and \; \mathcal{Y}_2=\mathcal{D}_1^2\mathcal{A} \nonumber
			\end{eqnarray}
			with $\mathcal{Y}_{i,j}^{(n)}=\left[\left(\mathcal{Y}_1\right)^{(n)}_{i,j},\, \left(\mathcal{Y}_2\right)^{(n)}_{i,j}\right]$ for $n\in \{1,2,...,n_3\}$, $i\in \{1,2,...,n_1\}$ and $j\in \{1,2,...,n_2\}$, ${\mathcal{D}_1^1}\mathcal{X}=\left[D_1^1\mathcal{X}^{(1)}|\,D_1^1\mathcal{X}^{(2)}|\,...| D_1^1\mathcal{X}^{(n_3)}|\, \right]$ and $\mathcal{D}_1^2\mathcal{X}=\left[D_1^2\mathcal{X}^{(1)}|\,D_1^2\mathcal{X}^{(2)}|\,...| D_1^2\mathcal{X}^{(n_3)}|\, \right]$.\\
			
			\noindent The  constrained optimization problem  \eqref{eq 3.4} can be written  as
			
			\begin{eqnarray}	\label{eq 3.5}
			&\underset{\mathcal{A},\, \mathcal{Z},\, \mathcal{W}}{min}& \;\left [ F\left(\mathcal{Z}\right)+G\left(\mathcal{W}\right) \right ]  \\
			& s.t & \mathcal{P}_\Omega (\mathcal{A})=\mathcal{P}_\Omega(\mathcal{M}),\; \mathcal{Z}=\mathcal{A},\; \mathcal{W}=\mathcal{D}_1\mathcal{A}  \nonumber
			\end{eqnarray}
			where 
			$$F\left(\mathcal{Z}\right)=\left\Vert \mathcal{Z}\right\Vert_*,\; G\left(\mathcal{W}\right)=\lambda\displaystyle{\sum_{n=1}^{n_3}\sum_{j=1}^{n_2}\sum_{i=1}^{n_1}\left\Vert \mathcal{Y}_{i,j}^{(n)}\right\Vert_2} , \mathcal{D}_1=\begin{pmatrix}
			\mathcal{D}_1^1 \\
			\mathcal{D}_1^2
			\end{pmatrix}\; and \; \mathcal{W}=\begin{pmatrix}
			\mathcal{W}_1\\
			\mathcal{W}_2
			\end{pmatrix}=\begin{pmatrix}
			\mathcal{Y}_1\\
			\mathcal{Y}_2
			\end{pmatrix}\in \mathbb{R}^{2n_1 \times n_2 \times n_3}.$$
			
			To solve the regularized optimisation problem \eqref{eq 3.5}, we can use the well known ADM method \cite{ hestenes1969multiplier,powell1969method}. The augmented Lagrangian associated to the problem  \eqref{eq 3.5} is given by
			\begin{eqnarray} \label{eq 3.6}
			L\left(\mathcal{A},\mathcal{Z},\mathcal{W},\mathcal{Q},\mathcal{B}\right) = F\left(\mathcal{Z}\right) + G\left(\mathcal{W}\right) &+& \left<\mathcal{A}-\mathcal{Z},\mathcal{Q}\right> + \dfrac{\beta_1}{2}\left\Vert \mathcal{A}-\mathcal{Z}\right\Vert_F^2 \nonumber \\
			&+& \left<\mathcal{D}_1\mathcal{A}-\mathcal{W},\mathcal{B}\right>+\dfrac{\beta_2}{2}\left\Vert \mathcal{D}_1\mathcal{A}-\mathcal{W} \right\Vert_F^2  
			\end{eqnarray}
			with $\mathcal{Q}\in \mathbb{R}^{n_1\times n_2 \times n_3}$, $\mathcal{B}=\begin{pmatrix}
			\mathcal{B}_1\\
			\mathcal{B}_2
			\end{pmatrix}\in \mathbb{R}^{2n_1\times n_2 \times n_3}$ are the Lagrangian multipliers and $\beta_1,\beta_2>0$ are the penalty parameters. 
			Therefore ADM leads to  the following sub-problems
			\begin{eqnarray}
			&\left(\mathcal{A}^k,\mathcal{Z}^k,\mathcal{W}^k\right)&=\underset{\mathcal{X},\mathcal{Z},\mathcal{W}}{arg\, min}\, L(\mathcal{A},\mathcal{Z},\mathcal{W},\mathcal{Q}^k,\mathcal{B}^k), \label{eq 3.7}\\
			&\mathcal{Q}^{k+1}&=\mathcal{Q}^k+\beta_1(\mathcal{A}^k-\mathcal{Z}^k),  \label{eq 3.8}\\
			&\mathcal{B}^{k+1}&=\mathcal{B}^k+\beta_2\left(\mathcal{D}_1\mathcal{A}^k-\mathcal{W}^k\right) \label{eq 3.9}.
			\end{eqnarray}
			
			\medskip
			\noindent Let us see now how to solve each of those sub-problems.\\
			\begin{itemize}
				\item \textbf{Solving the $\mathcal{A}$-problem :} From \eqref{eq 3.7} for a given $\mathcal{Z},\; \mathcal{W}$, we compute  the approximation $	\mathcal{A}^{k}$ by solving for $	\mathcal{A}$ the minimization problem 
				\begin{eqnarray*}
					\mathcal{A}^{k}=\underset{\mathcal{A}}{\arg\, min}\; \dfrac{\beta_1}{2}\left\Vert \mathcal{A}-\mathcal{Z}+\dfrac{\mathcal{Q}^k}{\beta_1}\right\Vert_F^2 + \dfrac{\beta_2}{2}\left\Vert \mathcal{D}_1\mathcal{A}-\mathcal{W} + \dfrac{\left(\mathcal{B}^k\right)}{\beta_2}\right\Vert_F^2 \nonumber 
				\end{eqnarray*}
				Then, the optimal value $\mathcal{A}^k$ satisfies  the following equation
				\begin{equation}\label{mat1}
				\beta_1 \mathcal{A}^k + \beta_2 \left(\mathcal{D}_1^1\right)^T \mathcal{D}_1^1 \mathcal{A}^k + \beta_2 \left(\mathcal{D}_1^2\right)^T \mathcal{D}_1^2 \mathcal{A}^k=\mathcal{R}
				\end{equation}
				with $\mathcal{R}=\beta_1 \mathcal{Z} - \mathcal{Q}^{k}+ \beta_2\left( \mathcal{D}_1^1\right)^T \mathcal{W}_1 -\left( \mathcal{D}_1^1 \right)^T \mathcal{B}^k_1 + \beta_2\left( \mathcal{D}_1^2\right)^T \mathcal{W}_2 - \left(\mathcal{D}_1^2\right)^T \mathcal{B}^k_2.$\\
				
				\noindent To solve the tensor equation \eqref{mat1},  we can   transform it to  a matrix equation by considering  the   frontal slices of the tensor $\mathcal{A}^k$. Then,  for $n\in \{1,2,...,n_3\}$, the matrix   $\left(\mathcal{A}^{k}\right)^{(n)}$  satisfies   the following matrix equation
				\begin{equation}
				\beta_1 \left(\mathcal{A}^{k}\right)^{(n)} + \beta_2(D_1^1)^TD_1^1\left(\mathcal{A}^{k}\right)^{(n)}+ \beta_2 (D_1^2)^TD_1^2\left(\mathcal{A}^{k}\right)^{(n)} =\mathcal{R}^{(n)}
				\label{eq 47}
				\end{equation}
				which can be also written as
				\begin{equation}
				\beta_1 \left(\mathcal{A}^k\right)^{(n)}+ \beta_2 \left(\mathcal{A}^k\right)^{(n)} \left(C_{n_2}^1\right)^TC_{n_2}^1+ \beta_2 \left(C_{n_1}^2\right)^T C_{n_1}^2\left(\mathcal{A}^k\right)^{(n)}=\mathcal{R}^{(n)} \label{eq 477}
				\end{equation}
				Since $C_{n_1}^2$ and $C_{n_2}^1$ are circulant matrices, they are diagonalizable using the discrete Fourier transformation, i.e., there exist $\Lambda_1$ and $\Lambda_2$ diagonal matrices such that
				$$C_{n_1}^2=F^*_{n_1}\Lambda_1 F_{n_1}, \;\; C_{n_2}^1=F^*_{n_2}\Lambda_2 F_{n_2}$$
				where $F_{n_1}$ and $F_{n_2}$ are the respectively the matrices representing the discrete Fourier transformation of size $n_1 \times n_1$ and $n_2 \times n_2$, then 
				$$\left(C_{n_1}^2\right)^TC_{n_1}^2=F^*_{n_1}\Lambda_1^2 F_{n_1}, \;\; \left(C_{n_2}^1\right)^TC_{n_2}^1=F^*_{n_2}\Lambda_2^2 F_{n_2}$$		
				By referring to \cite{ji2016tensor}, we can rewrite \eqref{eq 477}, for each $n\in \{1,2,...,n_3\}$,  as 
				$$\left(F_{n_2}^*\otimes F_{n_1}^*\right)\left(\beta_1 I\otimes I +\beta_2 \Lambda_2^2 \otimes I + \beta_2 I \otimes \Lambda_1^2  \right)\left(F_{n_2}\otimes F_{n_1}\right)vect\left(\left(\mathcal{A}^k\right)^{(n)}\right)=vect\left(\mathcal{R}^{(n)}\right)$$
				and since $\beta_1 I\otimes I +\beta_2 \Lambda_2^2 \otimes I + \beta_2 I \otimes \Lambda_1^2$ is an invertible matrix we get that
				\begin{equation} \label{eq A1}
				vect\left(\left(\mathcal{A}^k\right)^{(n)}\right)=\left(F_{n_2}^*\otimes F_{n_1}^*\right)\left(\beta_1 I\otimes I +\beta_2 \Lambda_2^2 \otimes I + \beta_2 I \otimes \Lambda_1^2  \right)^{-1}\left(F_{n_2}\otimes F_{n_1}\right)vect(\mathcal{R}^{(n)}).
				\end{equation}
				As the parameters $\beta_1$ and $\beta_2$ are strictly positive numbers, this shows that  for each $n\in \{1,2,...,n_3\}$, the equation \eqref{eq 477} has a unique solution.
				
				\item \textbf{Solving  the $\mathcal{Z}$-problem :} Given $\mathcal{X}$ and  $\mathcal{W}$,  the value of $\mathcal{Z}^k$ satisfies the following optimization problem
				\begin{eqnarray}
				\mathcal{Z}^{k}&=&\underset{\mathcal{Z}}{arg\, min}\; F\left(\mathcal{Z}\right) + \dfrac{\beta_1}{2}\left\Vert \mathcal{Z}-\mathcal{A}-\dfrac{\mathcal{Q}^k}{\beta_1}\right\Vert_F^2.\nonumber \\
				& = & \underset{\mathcal{Z}}{arg\, min}\; \left\Vert \mathcal{Z}\right\Vert_* + \dfrac{\beta_1}{2}\left\Vert \mathcal{Z}-\mathcal{A}-\dfrac{\mathcal{Q}^k}{\beta_1}\right\Vert_F^2. \nonumber
				\end{eqnarray}
				Then, from the result of Theorem \ref{theo 8}, we get
				\begin{equation}
				\mathcal{Z}^{k}=\mathcal{D}_\tau \left(\mathcal{A}+\dfrac{\mathcal{Q}^k}{\beta_1}\right)
				\label{eq 3.11}
				\end{equation}
				with $\tau = \dfrac{1}{\beta_1}$.
				
				\item \textbf{Solving the $\mathcal{W}$-problem :} For a given $\mathcal{X}$ and $ \mathcal{Z}$,  $\mathcal{W}^k$ is obtained by solving the following  sub-problems:  for $n\in \{1,2,...,n_3\}$ the $n^{th}$ sub-problem is given by
				\begin{eqnarray}
				\left(\mathcal{W}^{k}\right)^{(n)}=\underset{\mathcal{W}^n}{arg\, min}\; \lambda \sum_{j=1}^{n_2}\sum_{i=1}^{n_1} \left\Vert \mathcal{Y}_{i,j}^n\right\Vert_2 &+& \dfrac{\beta_2}{2}\left\Vert \mathcal{Y}_1^{(n)} - D_1^1\mathcal{A}^{(n)} - \dfrac{\left(\mathcal{B}^k_1\right)^n}{\beta_2}\right\Vert_F^2 \nonumber\\
				& + & \dfrac{\beta_2}{2}\left\Vert \mathcal{Y}_2^{(n)} - D_1^2\mathcal{A}^{(n)} - \dfrac{\left(\mathcal{B}^k_2\right)^n}{\beta_2}\right\Vert_F^2. \nonumber		\;\;\;\;
				\end{eqnarray}
				which is equivalent to solve the following $\mathcal{Y}$-problem 
				\begin{eqnarray}
				\left(\mathcal{Y}^{k}\right)^{(n)}_{i,j}=\underset{\mathcal{Y}^n_{i,j}}{arg\, min}\; \lambda  \left\Vert \mathcal{Y}_{i,j}^{(n)}\right\Vert_2 &+& \dfrac{\beta_2}{2}\left[ \left(\mathcal{Y}_1^{(n)}\right)_{i,j} - \left(D_1^1\mathcal{A}^{(n)}\right)_{i,j} - \dfrac{\left(\mathcal{B}^k_1\right)_{i,j}^{(n)}}{\beta_2}\right]^2 \nonumber\\
				& + & \dfrac{\beta_2}{2}\left[ \left(\mathcal{Y}_2^{(n)}\right)_{i,j} - \left(D_1^2\mathcal{A}^{(n)}\right)_{i,j} - \dfrac{\left(\mathcal{B}^k_2\right)_{i,j}^{(n)}}{\beta_2}\right]^2.\;\;\;\;\;\nonumber
				\end{eqnarray}
				By using the $2$-D shrinkage formula, we will get, for $1\leq i \leq n_1$ and $1\leq j \leq n_2$, the following expression
				\begin{equation}
				\left(\mathcal{Y}^{k}\right)^{(n)}_{i,j}=\max\{\left\Vert \mathcal{S}_{i,j}^{(n)}\right\Vert_2-\dfrac{\lambda}{\beta_2},0\}\dfrac{\mathcal{S}_{i,j}^{(n)}}{\left\Vert \mathcal{S}_{i,j}^{(n)}\right\Vert_2}
				\label{eq 3.12}
				\end{equation}
				where $\mathcal{S}_{i,j}^{(n)}=\left\{\left(D_1^1\mathcal{A}^{(n)}\right)_{i,j} + \dfrac{\left(\mathcal{B}^k_1\right)_{i,j}^{(n)}}{\beta_2},\, \left(D_1^2\mathcal{A}^{(n)}\right)_{i,j} + \dfrac{\left(\mathcal{B}^k_2\right)_{i,j}^{(n)}}{\beta_2} \right\}$ and we set $0\left(\dfrac{0}{0}\right)=0.$
			\end{itemize}
			The tensor completion procedure using the tensor nuclear norm and first order total variation (TNN-TV$_1$) is summarized in Algorithm \ref{alg:6}. 
			\begin{algorithm}[H]
				\caption{Tensor completion using the  tensor nuclear norm and the first order total variation (TNN-TV$_1$).}
				\label{alg:6}
				\begin{algorithmic}[1]
					\STATE \textbf{Initialize} $\mathcal{Z},\; \mathcal{W},\;\mathcal{Q},\;  \mathcal{B}$, $\lambda$, $\beta_1$ and $\beta_2$.
					\WHILE{not converged}
					\STATE  Update $\mathcal{A}^{k}$ from \eqref{eq A1} ;  
					\STATE Update $\mathcal{Z}^{k}$ from \eqref{eq 3.11};
					\STATE Update $\mathcal{W}^{k}$ from \eqref{eq 3.12};
					\STATE Update $\mathcal{Q}^{k+1}$ from \eqref{eq 3.8};
					\STATE Update $\mathcal{B}^{k+1}$ from\eqref{eq 3.9};
					\ENDWHILE
				\end{algorithmic}
			\end{algorithm}
		\end{subsection}
		
		\begin{subsection}{Tensor completion using tensor nuclear norm and the second order total variation}
			In this part we  apply the total variation with the second order derivative. Then  second proposed model is  formulated as 
			\begin{eqnarray}
			&\underset{\mathcal{A}}{min}& \, \left\Vert \mathcal{A}\right\Vert_* + \lambda \textbf{TV}_2(\mathcal{A})\nonumber\\
			&s.t& \mathcal{P}_\Omega (\mathcal{A})=\mathcal{P}_\Omega (\mathcal{M})
			\label{eq 3.13}
			\end{eqnarray}
			where $\mathcal{A}$, $\mathcal{M}$ $\left\Vert . \right\Vert_*$ and $\Omega$ play the same role as in the first model \eqref{eq 3.3}.  The expression of   $\textbf{TV}_2$ is given by 
			\begin{equation}
			\textbf{TV}_2\left(\mathcal{A}\right)=\left[TV_2(\mathcal{A}^{(1)})|\, TV_2(\mathcal{A}^{(2)})| ... | TV_2(\mathcal{A}^{(n_3)})\right] \in \mathbb{R}^{n_1 \times n_2 \times n_3} \nonumber
			\end{equation}
			where
			\begin{equation}
			TV_2\left(\mathcal{A}^{(n)}\right)=\sum_{i=1}^{n_1}\sum_{j=1}^{n_2}\sqrt{\left(D_2^1\mathcal{A}^{(n)}\right)^2_{i,j}+\left(D_2^2\mathcal{A}^{(n)}\right)^2_{i,j}} \nonumber
			\end{equation}
			for each $n\in \{1,2,...,n_3\}$. 
			The matrices  $D_2^1$ and $D_2^2$ are the second derivative operators in the first and in the second direction, respectively , satisfying  for $n\in \{1,2,...,n_3\}$
			\begin{equation}
			D_2^1\mathcal{A}^{(n)}=\mathcal{A}^{(n)}C_{n_2} \;\; and \;\; D_2^2\mathcal{A}^{(n)}=C_{n_1}\mathcal{A}^{(n)} \nonumber
			\end{equation}
			with 
			$$C_i=\dfrac{1}{2}\begin{pmatrix}
			-2 & 1 & 0 & 0 & ... & 1 \\
			1 & -2 & 1 & 0 & ... & 0 \\
			0 & \ddots &  \ddots & \ddots & \ddots & \vdots \\
			\vdots & \ddots &  \ddots & \ddots & \ddots & \vdots\\
			0 & ... &... & 1 & -2 & 1\\
			1 & ... &... & 0 & 1 & -2 
			\end{pmatrix}\in \mathbb{R}^{i\times i};\; i=n_1,\,  n_2.$$
			The optimization problem  \eqref{eq 3.13} is equivalent to the following one
			\begin{eqnarray}
			&\underset{\mathcal{A},\, \mathcal{Z},\, \mathcal{Y}}{min}& \;\left [ \left\Vert\mathcal{Z}\right\Vert_*+\lambda\sum_{n=1}^{n_3}\sum_{j=1}^{n_2}\sum_{i=1}^{n_1}\left\Vert \mathcal{Y}_{i,j}^{(n)}\right\Vert_2 \right] \nonumber \\
			&s.t& \mathcal{P}_\Omega \left(\mathcal{A}\right)=\mathcal{P}_\Omega(\mathcal{M}),\; \mathcal{Z}=\mathcal{A},\; \mathcal{Y}_1=\mathcal{D}_2^1 \mathcal{A} \; and \; \mathcal{Y}_2=\mathcal{D}_2^2 \mathcal{A}
			\label{eq 3.14}
			\end{eqnarray}
			with $\mathcal{Y}_{i,j}^{(n)}=\left[ \left(\mathcal{Y}_1\right)_{i,j}^n,\,  \left( \mathcal{Y}_2 \right)_{i,j}^{(n)} \right]$ for $n\in \{1,2,...,n_3\}$, $i\in \{1,2,...,n_1\}$ and $j\in \{1,2,...,n_2\}$, $\mathcal{D}_2^1\mathcal{A}=\left[D_2^1 \mathcal{A}^{(1)}|\, D_2^1 \mathcal{A}^{(2)}|\, .... | D_2^1 \mathcal{A}^{(n_3)}\right]$ and $\mathcal{D}_2^2\mathcal{A}=\left[D_2^2 \mathcal{A}^{(1)}|\, D_2^2 \mathcal{A}^{(2)}|\, .... | D_2^2 \mathcal{A}^{(n_3)}\right]$. \\
			
			\noindent The  constrained optimization problem \eqref{eq 3.14}  is transformed to the following one
			\begin{eqnarray}
			&\underset{\mathcal{A},\mathcal{Z},\mathcal{W}}{arg\, min}& \, \left[F(\mathcal{Z})+G(\mathcal{W})\right] \nonumber \\
			& s.t & \mathcal{P}_{\Omega} (\mathcal{A})=\mathcal{P}_{\Omega}(\mathcal{M}),\; \mathcal{A}=\mathcal{Z},\; \mathcal{D}_2\mathcal{A}=\mathcal{W} \label{eq 3.15}
			\end{eqnarray}
			where $F\left(\mathcal{Z}\right)=\left\Vert \mathcal{Z}\right\Vert_*$, $ G\left(\mathcal{W}\right)=\lambda \displaystyle{\sum_{n=1}^{n_3}\sum_{j=1}^{n_2}\sum_{i=1}^{n_1}\left\Vert \mathcal{Y}_{i,j}^{(n)}\right\Vert_2}$, $\mathcal{D}_2=\begin{pmatrix}
			\mathcal{D}_2^1\\
			\mathcal{D}_2^2
			\end{pmatrix}$ and $\mathcal{W}=\begin{pmatrix}
			\mathcal{W}_1\\
			\mathcal{W}_2
			\end{pmatrix}=\begin{pmatrix}
			\mathcal{Y}_1\\
			\mathcal{Y}_2
			\end{pmatrix}$. \\
			We  observe that \eqref{eq 3.15} is similar to \eqref{eq 3.5}. Thus, we can use the same procedure of solving \eqref{eq 3.5} to solve \eqref{eq 3.15}.
			The augmented Lagrangian associated to the  optimization problem \eqref{eq 3.15} is given by
			\begin{eqnarray}
			L(\mathcal{A},\mathcal{Z},W,\mathcal{Q},B)=F(\mathcal{Z}) + G(\mathcal{W}) &+& \left< \mathcal{A}-\mathcal{Z},\mathcal{Q}\right> + \dfrac{\beta_1}{2}\left\Vert \mathcal{A}-\mathcal{Z}\right\Vert_F^2 \nonumber\\
			&+&\left<\mathcal{D}_2\mathcal{A}-\mathcal{W},\mathcal{B}\right>+\dfrac{\beta_2}{2}\left\Vert \mathcal{D}_2\mathcal{A}-\mathcal{W}\right\Vert_F^2. \label{3.16}
			\end{eqnarray}
			Therefore, using ADM, we have to solve  the following sub-problems
			\begin{eqnarray}
			&\left(\mathcal{A}^k,\mathcal{Z}^k,\mathcal{W}^k\right)&=\underset{\mathcal{A},\mathcal{Z},\mathcal{W}}{arg\, min}\, L(\mathcal{A},\mathcal{Z},\mathcal{W},\mathcal{Q}^k,\mathcal{B}) \label{eq 3.17} \\
			&\mathcal{Q}^{k+1}=&\mathcal{Q}^k + \beta_1 \left(\mathcal{A}^{k}-\mathcal{Z}^{k}\right) \label{eq 3.18}\\
			&\mathcal{B}^{k+1}=&\mathcal{B}^k + \beta_2 \left(\mathcal{D}_2\mathcal{A}^{k}-\mathcal{W}^{k}\right).  \label{eq 3.19}
			\end{eqnarray}
			Let us see now how to solve each of those sub-problems.
			\begin{itemize}
				\item \textbf{The $\mathcal{A}$-problem: }For fixed $\mathcal{Z}$ and $\mathcal{W}$, each  frontal slice of the approximation $\mathcal{A}^k$ satisfies   the Sylvester matrix equation
				\begin{equation}
				\beta_1 \left(\mathcal{A}^{k}\right)^{(n)} + \beta_2(D_2^1)^TD_2^1\left(\mathcal{A}^{k}\right)^{(n)}+ \beta_2 (D_2^2)^TD_2^2\left(\mathcal{A}^{k}\right)^{(n)} =\mathcal{R}^{(n)}
				\label{eq A2}
				\end{equation}			
				where  $\mathcal{R}=\beta_1 \mathcal{Z}-\mathcal{Q}^k+\beta_2 \left(\mathcal{D}_2^1\right)^T \mathcal{W}_1-\left(\mathcal{D}_2^1\right)^T \mathcal{B}_1^k+\beta_2 \left(\mathcal{D}_2^2\right)^T\mathcal{W}_2-\left(\mathcal{D}_2^2\right)^T\mathcal{B}_2^k$,\\
				which can be written  as 
				$$\beta_1 \left(\mathcal{A}^k\right)^{(n)}+ \beta_2 \left(\mathcal{A}^{k}\right)^{(n)}C_{n_2}^T C_{n_2}+\beta_2 C_{n_1}^TC_{n_1}\left(\mathcal{A}^{k}\right)^{(n)}=\mathcal{R}^{(n)}. $$
				Using the same idea as for  \ref{eq 477},  each frontal slice of $\mathcal{A}^k$ satisfies 
				\begin{equation}\label{eq A3}
				vect\left(\left(\mathcal{A}^k\right)^{(n)}\right)=\left(F_{n_2}^*\otimes F_{n_1}^*\right)\left(\beta_1 I\otimes I +\beta_2 \Lambda_2^2 \otimes I + \beta_2 I \otimes \Lambda_1^2  \right)^{-1}\left(F_{n_2}\otimes F_{n_1}\right)vect(\mathcal{R}^{(n)}).
				\end{equation} 
				with 
				$$\Lambda_1= F_{n_1}C_{n_1} F_{n_1}^*,\; {\rm and} \;  \Lambda_2= F_{n_2}C_{n_2} F_{n_2}^*,$$
				where $F_{n_i}$ is the Fourier matrix of size $n_i \times n_i$  for $i=1,2$.
				
				\item \textbf{The $\mathcal{Z}$-problem: }For $\tau=\dfrac{1}{\beta_1}$ and for a given $\mathcal{A}$ and $\mathcal{W}$ we get
				\begin{equation} \label{eq 3.21}
				\mathcal{Z}^{k}=\mathcal{D}_\tau\left(\mathcal{A}+\dfrac{\mathcal{Q}^k}{\beta_1}\right).
				\end{equation} 
				\item \textbf{The $\mathcal{W}$-problem: }By applying 2D shrinkage formula on each frontal slice of $\mathcal{Y}$ for a given $\mathcal{A}$ and $\mathcal{Z}$,  we  get  
				\begin{equation} \label{eq 3.22}
				\left(\mathcal{Y}^{k}\right)^{(n)}_{i,j}=max \left\{\left\Vert \mathcal{S}^{(n)}_{i,j}\right\Vert_2 - \dfrac{\lambda}{\beta_2},0\right\}\dfrac{\mathcal{S}^{(n)}_{i,j}}{\left\Vert \mathcal{S}^{(n)}_{i,j}\right\Vert_2}
				\end{equation}
				where $\mathcal{S}^{(n)}_{i,j}=\left\{\left(D_2^1\mathcal{A}^{(n)}\right)_{i,j}+\dfrac{\left(\mathcal{B}^k_1\right)^{(n)}_{i,j}}{\beta_2}, \left(D_2^2\mathcal{A}^{(n)}\right)_{i,j}+\dfrac{\left(\mathcal{B}^k_2\right)^{(n)}_{i,j}}{\beta_2}\right\}$ with  $0\left(\dfrac{0}{0}\right)=0.$	
			\end{itemize}
			The different steps of the 	tensor completion using the tensor nuclear norm and total variation (TNN-TV$_2$) is summarized in the following algorithm. 
			\begin{algorithm}[H]
				\caption{Tensor completion using the tensor nuclear norm and the second order total variation (TNN-TV$_2$).}
				\label{alg:7}
				\begin{algorithmic}[1]
					\STATE \textbf{Initialize} $\mathcal{Z},\;\mathcal{W},\;\mathcal{Q},\;  \mathcal{B}$, $\lambda$, $\beta_1$ and $\beta_2$.
					\WHILE{not converged}
					\STATE Update $\mathcal{A}^{k}$ from \eqref{eq A3} ;  
					\STATE Update $\mathcal{Z}^{k}$ from \eqref{eq 3.21};
					\STATE Update $\mathcal{W}^{k}$ from \eqref{eq 3.22};
					\STATE Update $\mathcal{Q}^{k+1}$ from \eqref{eq 3.18};
					\STATE Update $\mathcal{B}^{k+1}$ from \eqref{eq 3.19};
					\ENDWHILE
				\end{algorithmic}
			\end{algorithm}
		\end{subsection}
		Now we discuss the complexity of the two algorithms TNN-TV$_1$ and TNN-TV$_2$ . As we used the fast Fourier transform, the cost of computing the $\mathcal{A}$-sub-problem  is $O(n_3n_1n_2 log(n_1n_2))$. The cost of computing $\mathcal{Z}$ in \eqref{eq 3.11} and \eqref{eq 3.21} is $O(\dfrac{n_3}{2}(2n_1^2n_2 +n_1n_2^2) )$. Computing  $\mathcal{W}$ in \eqref{eq 3.12} and \eqref{eq 3.22} requires  $O(n_1n_2n_3)$ arithmetic  operations. 
	\end{section}
	
	\section{Convergence analysis}
	In this section we study the convergence of the proposed approaches. As the two methods are similar, we will give theoretical results only for sequences obtained by Algorithm \ref{alg:6}. 
	Notice first  that the functions $F$ and $G$ defined earlier are closed, proper and convex. Then, thanks to \cite{facchinei2007finite}\cite{rockafellar1970convex},  the optimization problem  \eqref{eq 3.5} is solvable, i.e., there exist $\mathcal{Z}^*$ and $\mathcal{W}^*$ not necessarily unique that minimize \eqref{eq 3.5}.\\
	Let us define the space $\mathbb{E}=\mathbb{R}^{n_1\times  n_2 \times n_3}\times \mathbb{R}^{n_1\times n_2 \times n_3}\times \mathbb{R}^{2n_1\times \times n_2 \times n_3}\times \mathbb{R}^{n_1\times \times n_2 \times n_3}\times \mathbb{R}^{2n_1\times \times n_2 \times n_3}$ which is closed and nonempty. We first recall the following theorem. 
	\begin{theo} \cite{wu2011augmented} 
		$\mathcal{A}^*$ is a solution of \eqref{eq 3.3}    if and only if there exist $\left(\mathcal{Z}^*,\mathcal{W}^*\right)\in \mathbb{R}^{n_1 \times n_2 \times n_3}\times \mathbb{R}^{2N_1 \times n_2 \times n_3}$ and $\left(\mathcal{Q}^*,\mathcal{B}^*\right))\in \mathbb{R}^{n_1 \times n_2 \times n_3}\times \mathbb{R}^{2N_1 \times n_2 \times n_3}$ such that $\left(\mathcal{A}^*,\mathcal{Z}^*,\mathcal{W}^*,\mathcal{Q}^*,\mathcal{B}^*\right)\in \mathbb{E}$ is a saddle point of $L$
		,i.e. 
		\begin{eqnarray} \label{eq 4.1}
		&L(\mathcal{A}^*,\mathcal{Z}^*,\mathcal{W}^*,\mathcal{Q},\mathcal{B})&\leq 	L(\mathcal{A}^*,\mathcal{Z}^*,\mathcal{W}^*,\mathcal{Q}^*,\mathcal{B}^*)\leq L(\mathcal{A},\mathcal{Z},\mathcal{W},\mathcal{Q}^*,\mathcal{B}^*) \nonumber  \\
		& \forall \left(\mathcal{A},\mathcal{Z},\mathcal{W},\mathcal{Q},\mathcal{B}\right)\in \mathbb{E}.&
		\end{eqnarray}
	\end{theo}
	The next theorem gives some convergence results on the sequences obtained from Algorithm \ref{alg:6}. 
	\begin{theo}
		Assume that $(\mathcal{A}^*,\mathcal{Z}^*,\mathcal{W}^*,\mathcal{Q}^*,\mathcal{B}^*)$ is a saddle point of $L$. The sequence $(\mathcal{X}^k,\mathcal{Z}^k,\mathcal{W}^k,\mathcal{Q}^k,\mathcal{B}^k)$ generated by Algorithm \ref{alg:6} satisfies:
		\begin{enumerate}
			\item $\underset{k \to +\infty}{lim} F(\mathcal{Z}^k)+G(\mathcal{W}^k)=F(\mathcal{Z}^*)+G(\mathcal{W}^*).$
			\item $\underset{k \to +\infty}{lim} \left\Vert \mathcal{A}^k - \mathcal{Z}^k \right\Vert=0.$
			\item $\underset{k \to +\infty}{lim} \left\Vert \mathcal{D}_1 \mathcal{A}^k - \mathcal{W}^k \right\Vert=0.$
		\end{enumerate}
	\end{theo}
	\begin{pf}
		From the first inequality of \eqref{eq 4.1} we get
		\begin{equation}
		\left<\mathcal{A}^*-\mathcal{Z}^*,\mathcal{Q}\right> + \left<\mathcal{D}_1\mathcal{A}^*-\mathcal{W}^*,\mathcal{B}\right> \leq \left<\mathcal{A}^*-\mathcal{Z}^*,\mathcal{Q}^*\right> + \left<\mathcal{D}_1\mathcal{A}^*-\mathcal{W}^*,\mathcal{B}^*\right> \;\; \forall (\mathcal{Q},\mathcal{W})
		\end{equation}
		which gives 
		$$\left\{
		\begin{array}{ll}
		\mathcal{A}^*=\mathcal{Z}^* \\
		\mathcal{D}_1\mathcal{A}^*=\mathcal{W}^*.
		\end{array}
		\right.$$\\
		Let us define the following quantities \\
		$$\bar{\mathcal{X}}^k=\mathcal{A}^k-\mathcal{A}^*,\; \bar{\mathcal{Z}}^k=\mathcal{Z}^k-\mathcal{Z}^*,\; \bar{\mathcal{W}}^k=\mathcal{W}^k-\mathcal{W}^*,\; \bar{\mathcal{Q}}^k=\mathcal{Q}^k-\mathcal{Q}^*,\; \bar{\mathcal{B}}^k=\mathcal{B}^k-\mathcal{B}^*.$$
		The main idea of the  proof is to show that the sequence $\left( \beta_2 \left\Vert \bar{\mathcal{Q}}^k \right\Vert_F^2 + \beta_1\left\Vert \bar{\mathcal{B}}^k \right\Vert_F^2  \right)_{k\geq 0}$ is decreasing. Notice that\\
		$$\bar{\mathcal{Q}}^{k+1}=\bar{\mathcal{Q}}^k-\beta_1\left(\bar{\mathcal{A}}^k-\bar{\mathcal{Y}}^k \right),\; \bar{\mathcal{W}}^{k+1}=\bar{\mathcal{W}}^k-\beta_2\left(\mathcal{D}_1\bar{\mathcal{A}}^k-\bar{\mathcal{B}}^k \right).$$
		We have 
		\begin{eqnarray}
		\left( \left\Vert \bar{\mathcal{Q}}^k \right\Vert_F^2 + \left\Vert \bar{\mathcal{B}}^k \right\Vert_F^2  \right) - \left( \left\Vert \bar{\mathcal{Q}}^{k+1} \right\Vert_F^2 + \left\Vert \bar{\mathcal{B}}^{k+1} \right\Vert_F^2  \right)&=& - 2 \beta_1 \left<\bar{\mathcal{Q}}^k, \bar{\mathcal{A}}^k - \bar{\mathcal{Y}}^k\right> - 2 \beta_1 \left<\bar{\mathcal{B}}^k, \mathcal{D}_1\bar{\mathcal{A}}^k - \bar{\mathcal{W}}^k\right> \nonumber \\
		& - & \beta_1^2  \left\Vert \bar{\mathcal{A}}^k -\bar{\mathcal{Y}}^k \right\Vert_F^2 - \beta_2^2  \left\Vert \mathcal{D}_1\bar{\mathcal{A}}^k -\bar{\mathcal{W}}^k \right\Vert_F^2.
		\end{eqnarray}
		From the second inequality of \eqref{eq 4.1} we obtain  the following inequalities for $(\mathcal{A},\mathcal{Z},\mathcal{W})=(\mathcal{A}^k,\mathcal{Z}^k,\mathcal{W}^k)$
		\begin{equation*}
		\;\;\;\;\;\;\; \left<\mathcal{A}^k-\mathcal{Z}^*,\mathcal{Q}^*\right>+\left<\mathcal{D}_1\mathcal{A}^k-\mathcal{W}^*,\mathcal{B}^*\right> + \beta_1\left< \mathcal{A}^k-\mathcal{Z}^*,\mathcal{A}^k-\mathcal{A}^* \right> + \beta_2\left< \mathcal{D}_1\mathcal{A}^k - \mathcal{W}^*,\mathcal{D}_1(\mathcal{A}^k-\mathcal{A}^*) \right> \geq 0
		\end{equation*}
		\begin{equation*}
		F(\mathcal{Z}^k)-F(\mathcal{Z}^*) + \left<\mathcal{A}^*-\mathcal{Y}^k,\mathcal{Q}^*\right>+\beta_1\left< \mathcal{A}^*-\mathcal{Y}^k,\mathcal{Y}^*-\mathcal{Y}^k \right> \geq 0.
		\end{equation*}
		\begin{equation*}
		G(\mathcal{W}^k)-F(\mathcal{W}^*) + \left<\mathcal{D}_1\mathcal{A}^*-\mathcal{W}^k,\mathcal{B}^*\right>+\beta_2\left< \mathcal{D}_1\mathcal{A}^*-\mathcal{W}^k,\mathcal{W}^*-\mathcal{W}^k \right> \geq 0.
		\end{equation*}
		Using \eqref{eq 3.7}, we get for $(\mathcal{A},\mathcal{Z},\mathcal{W})=(\mathcal{A}^k,\mathcal{Z}^k,\mathcal{W}^k)$ that
		\begin{equation*}
		\left<\mathcal{A}^*-\mathcal{Z}^k,\mathcal{Q}^k\right>+\left<\mathcal{D}_1\mathcal{A}^*-\mathcal{W}^k,\mathcal{B}^k\right> + \dfrac{\beta_1}{2}\left< \mathcal{A}^*-\mathcal{Z}^k,\mathcal{Z}^*-\mathcal{Z}^k \right> + \beta_2\left< \mathcal{D}_1\mathcal{A}^* - \mathcal{W}^k,\mathcal{W}^*-\mathcal{W}^k \right> \geq 0
		\end{equation*}
		\begin{equation*}
		F(\mathcal{Z}^*)-F(\mathcal{Z}^k) + \left<\mathcal{A}^k-\mathcal{Y}^*,\mathcal{Q}^k\right>+\beta_1\left< \mathcal{X}^k-\mathcal{Y}^*,\mathcal{A}^k-\mathcal{A}^* \right>\geq 0.
		\end{equation*}
		\begin{equation*}
		G(\mathcal{W}^*)-G(\mathcal{W}^k) + \left<\mathcal{D}_1\mathcal{A}^k-\mathcal{W}^*,\mathcal{B}^k\right>+\beta_2\left< \mathcal{D}_1\mathcal{A}^k-\mathcal{W}^*,\mathcal{D}_1(\mathcal{A}^k-\mathcal{A}^*) \right>\geq 0.
		\end{equation*}
		By regrouping terms we  get
		\begin{equation*}
		-  \left<\bar{\mathcal{Q}}^k, \bar{\mathcal{A}}^k - \bar{\mathcal{Y}}^k\right> -   \left<\bar{\mathcal{B}}^k, \mathcal{D}_1\bar{\mathcal{A}}^k - \bar{\mathcal{W}}^k\right> \geq  
		\beta_1  \left\Vert \bar{\mathcal{A}}^k -\bar{\mathcal{Y}}^k \right\Vert_F^2 + \beta_2  \left\Vert \mathcal{D}_1\bar{\mathcal{A}}^k -\bar{\mathcal{W}}^k \right\Vert_F^2
		\end{equation*}
		then
		\begin{equation*}
		\left( \beta_2 \left\Vert  \bar{\mathcal{Q}}^k \right\Vert_F^2 + \beta_1 \left\Vert \bar{\mathcal{B}}^k \right\Vert_F^2  \right) - \left(\beta_2 \left\Vert \bar{\mathcal{Q}}^{k+1} \right\Vert_F^2 + \beta_1\left\Vert \bar{\mathcal{B}}^{k+1} \right\Vert_F^2  \right)\geq
		\beta_1^2\beta_2  \left\Vert \bar{\mathcal{A}}^k -\bar{\mathcal{Y}}^k \right\Vert_F^2 + \beta_2^2\beta_1  \left\Vert \mathcal{D}_1\bar{\mathcal{A}}^k -\bar{\mathcal{W}}^k \right\Vert_F^2.
		\end{equation*}
		Thus, the sequence $\left(\beta_2\left\Vert \bar{\mathcal{Q}}^k \right\Vert_F^2 + \beta_1 \left\Vert \bar{\mathcal{B}}^k \right\Vert_F^2\right)_{k\geq 0}$ is decreasing, which gives \\
		\begin{equation*}
		\sum_{k= 0}^{+\infty}\left(\beta_1^2\beta_2  \left\Vert \bar{\mathcal{A}}^k -\bar{\mathcal{Y}}^k \right\Vert_F^2 + \beta_2^2\beta_1  \left\Vert \mathcal{D}_1\bar{\mathcal{A}}^k -\bar{\mathcal{W}}^k \right\Vert_F^2\right) \leq \beta_2 \left\Vert \bar{\mathcal{Q}}^0 \right\Vert_F^2 +\beta_1 \left\Vert \bar{\mathcal{W}}^0\right\Vert_F^2.
		\end{equation*}
		Therefore\\
		$$\left\{
		\begin{array}{lll}
		\left(\mathcal{Q}^k\right)_{k\geq 0}\; and\; \left(\mathcal{B}^k\right)_{k\geq 0}\; are \; bounded \\
		\underset{k \to +\infty}{lim} \left\Vert \mathcal{A}^k - \mathcal{Z}^k \right\Vert=0.\\
		\underset{k \to +\infty}{lim} \left\Vert \mathcal{D}_1 \mathcal{A}^k - \mathcal{W}^k \right\Vert=0.
		\end{array}
		\right.$$\\
		In addition, by using again the second inequality of \eqref{eq 4.1} for $(\mathcal{A},\mathcal{Z},\mathcal{W})=(\mathcal{A}^k,\mathcal{Z}^k,\mathcal{W}^k)$ we obtain 
		\begin{eqnarray*}
			F(\mathcal{Z}^*)+G(\mathcal{W}^*)\leq F(\mathcal{Z}^k)+G(\mathcal{W}^k) &+& \left< \mathcal{A}^k-\mathcal{Z}^k,\mathcal{Q}^k \right> + \dfrac{\beta_1}{2}\left\Vert \mathcal{A}^k-\mathcal{Z}^k\right\Vert_F^2 \nonumber \\
			&+& \left<\mathcal{D}_1\mathcal{A}^k-\mathcal{W}^k,\mathcal{B}^k\right>+\dfrac{\beta_2}{2}\left\Vert \mathcal{D}_1\mathcal{A}^k-\mathcal{W}^k \right\Vert_F^2,
		\end{eqnarray*}
		and 
		\begin{eqnarray*}
			F(\mathcal{Z}^*)+G(\mathcal{W}^*)\geq F(\mathcal{Z}^k)+G(\mathcal{W}^k)&+& \left<\mathcal{A}^k-\mathcal{Z}^k,\mathcal{Q}^k\right> + \dfrac{\beta_1}{2}\left\Vert \mathcal{A}^k-\mathcal{Z}^k\right\Vert_F^2 \nonumber \\
			& + & \left<\mathcal{D}_1\mathcal{A}^k-\mathcal{W}^k,\mathcal{B}^k\right>+\dfrac{\beta_2}{2}\left\Vert \mathcal{D}_1\mathcal{A}^k-\mathcal{W}^k\right\Vert_F^2.
		\end{eqnarray*}
		Hence 
		\begin{equation*}
		\underset{k\to +\infty}{lim\, sup}\;  F(\mathcal{Z}^k)+G(\mathcal{W}^k)\leq F(\mathcal{A}^*)+G(\mathcal{W}^*)\leq 	\underset{k\to +\infty}{lim\, inf} F(\mathcal{Z}^k)+G(\mathcal{W}^k).
		\end{equation*}
	\end{pf}
	\section{Numerical experiments}
	\label{sec 4}
	In this section, we  give some numerical tests to show the performance of our proposed algorithms TNN-TV$_1$ and TNN-TV$_2$  and  compare  them with the results obtained by other known   algorithms for image and video  completion,  such as  TNN \cite{zhang2014novel},  SiLRTC-TT \cite{bengua2017efficient} and MF-TV \cite{ji2016tensor}.\\
	The quality of the recovered images is measured by  computing the relative squared error (RSE), and the peak signal-to-noise-ration (PSNR), defined by
	\begin{equation*}
	RSE=\dfrac{\left\Vert \mathcal{A}_{ori}-\mathcal{A}\right\Vert_F^2}{\left\Vert \mathcal{A}\right\Vert_F^2},
	\end{equation*}
	and 
	\begin{equation*}
	PSNR=10 \,log_{10} \dfrac{Max_{\mathcal{A}}^2}{\left\Vert \mathcal{A}-\mathcal{A}_{ori}\right\Vert_F^2},
	\end{equation*}
	where $\mathcal{A}_{ori}$ is the original tensor, $\mathcal{A}$ is the recovered tensor and $Max_{\mathcal{A}} $ is the maximum pixel value of $\mathcal{A}$. The convergence  stopping criterion is defined by computing the relative error of $\mathcal{A}$ between two successive iterations as follows
	\begin{equation}
	\dfrac{\left\Vert \mathcal{A}^{k+1}-\mathcal{A}^k\right\Vert_F^2}{\left\Vert \mathcal{A}^k\right\Vert_F^2}\leq 10^{-4}.
	\end{equation} 
	In all the experiments we used fixed  values of regularization  and  penalty parameters. For the algorithm TNN-TV$_1$ we  used  $\lambda=0.1,\; \beta_1=0.01$ and $\beta_2=0.0001$, and for the algorithm TNN-TV$_2$ we set  $\lambda=1$, $\beta_1=0.01$ and $\beta_2=0.0001$.
	\begin{subsection}{Images}
		For this example,   we used color  images of size $256\times 256 \times 3$.  
		In Figure \ref{fig 1}, we  reported  the obtained visual results of  TNN,  SiLRTC-TT, MF-TV, TNN-TV$_1$ and TNN-TV$_2$,   with $SR=0.1$ , where $SR$ represents  the percentage of the data remained in the image. In Table \ref{tab 1} we compared   the efficiency of our two algorithms with TNN,   SiLRTC-TT and MF-TV by comparing the values of $RSE$ and $PSNR$. 
		\begin{figure}[H]
			\centering
			\begin{tabular}{c}
				\includegraphics[width=0.9\linewidth]{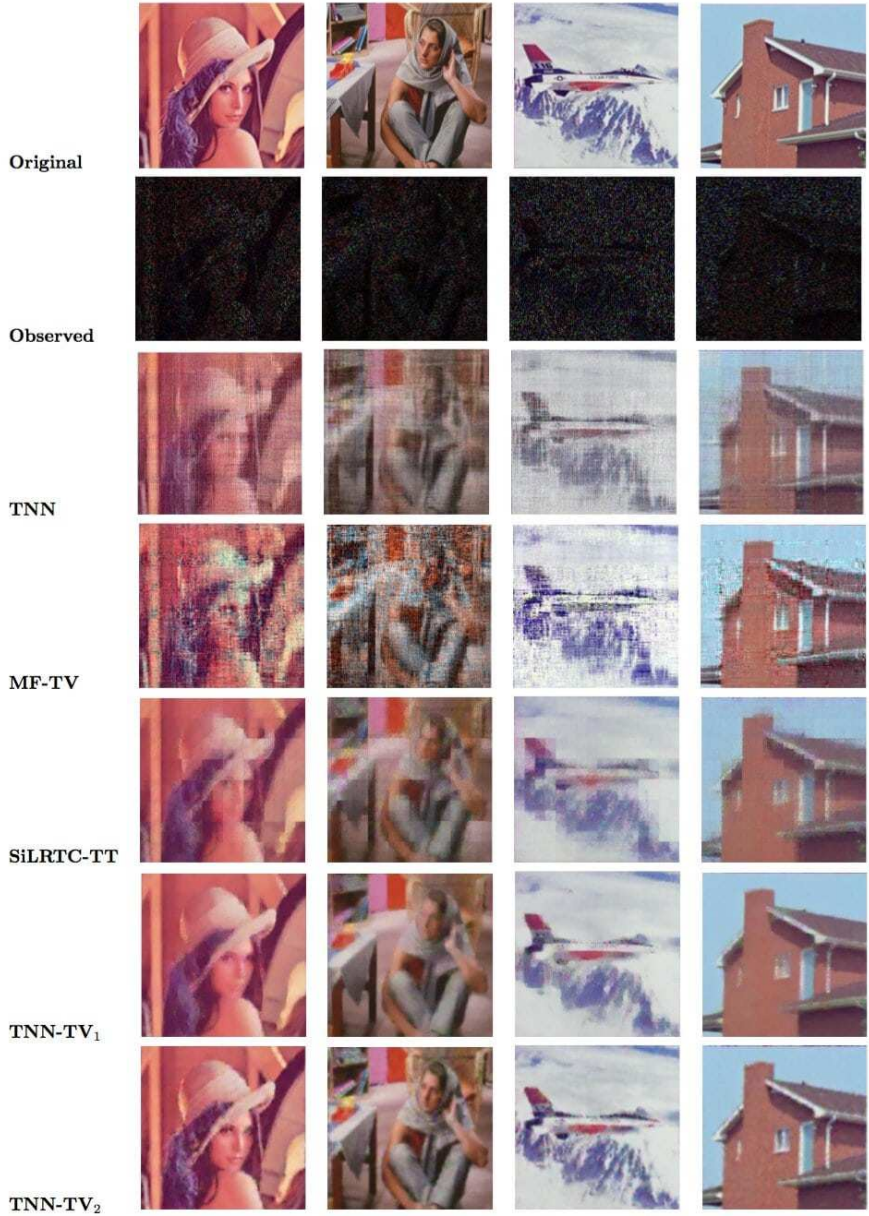}
			\end{tabular}
			\caption{The results of the algorithms TNN, MF-TV, ,SiLRTC-TT, TNN-TV$_1$ and TNN-TV$_2$  for the  images "Lena", "Barbara", "airplane" and "house" for $SR=0.1$.} 
			\label{fig 1}
		\end{figure}
	
		\noindent 	From Figure \ref{fig 1} and Table \ref{tab 1},  we can see the efficiency of our algorithms as compared  the others for different images.  We remark also that   TV$_2$ gives better results than TV$_1$. 
	
			\begin{table}
			\centering
			\footnotesize\begin{tabular}{l|l|l|l|l|l|l|l|l|l}
				\hline \multicolumn{2}{c|}{images} & 	\multicolumn{2}{|c|}{Lena} &  \multicolumn{2}{|c|}{Barbara} & \multicolumn{2}{|c|}{airplane} & \multicolumn{2}{|c}{house}   \\
				\hline $SR$ & algorithms & RSE & PSNR & RSE & PSNR & RSE & PSNR & RSE & PSNR \\
				\hline \multirow{5}{*}{0.1}  & TNN      &  0.1942   & 19.2585   &  0.2362   & 18.8644   &  0.1415  & 18.2747  & 0.1596  & 20.2903   \\
				&   SiLRTC-TT    &  0.1521   & 21.4342   &   0.2020  &  20.2838  &  0.1338  & 19.4385  & 0.1451     & 20.9854   \\
				&   MF-TV    &  0.3345   &  14.6340  & 0.3867    & 14.4866   &  0.1946  & 16.1082  &  0.1867   &  18.8811  \\
				&    TNN-TV$_1$   &  0.1203   & 23.5132   &  0.1560   &  22.3721  &  0.1102  & 21.0482 &  0.1080   &  23.6369  \\
				&   TNN-TV$_2$    & \textbf{0.0975}    &  \textbf{25.3402}  &  \textbf{0.1269}   & \textbf{24.1606}   &  \textbf{0.0869}  & \textbf{23.2522}  &  \textbf{0.0896}   &  \textbf{25.4389}  \\
				\hline \multirow{5}{*}{0.2}  &   TNN    &  0.1271       &  23.0104  &  0.2362        & 18.8644    &  0.1415  & 18.2747   & 0.1005    & 24.3883  \\
				&   SiLRTC-TT    &   0.1077  &  24.4344  & 0.2020    &  20.2838  &  0.1338  & 19.4385  &  0.1000  &  24.3342    \\
				&   MF-TV    &  0.1040   &  24.7759  &  0.3867   &  14.4866  & 0.1946   & 16.1082  &  0.0993      & 24.3628    \\
				&   TNN-TV$_1$    &  0.0837   & 26.6646   &  0.1560   & 22.3721   &  0.1102  & 21.0482  &  \textbf{0.0695}    & \textbf{27.5520}  \\
				&   TNN-TV$_2$    &    \textbf{0.0733}  &  \textbf{27.8160}  &  \textbf{0.1269}   & \textbf{24.1606}   &  \textbf{0.0869}  & \textbf{23.2522} &  0.0703   &   27.5507   \\
				\hline
			\end{tabular}
			\captionof{table}{The values of RSE and PSNR for  TNN,  SiLRTC-TT, MF-TV, TNN-TV$_1$ and TNN-TV$_2$ with the  images "Lena", "Barbara", "airplane" and "house" using $SR=0.1,\,0.2. $}
			\label{tab 1}
		\end{table}

		\begin{table}
			\centering
			\begin{tabular}{l|l|l|l|l|l|l|l|l}
				\hline SR & \multicolumn{4}{|c|}{0.2}& \multicolumn{4}{|c}{0.3} \\
				\hline RSE & 0.13 & 0.11 & 0.09 & 0.08 & 0.1 & 0.09 & 0.08 & 0.07 \\
				\hline
				\hline TNN &  3.0392      &     -    &  -     &  -      &    1.6736     &    -     &   -    &  -       \\
				\hline TNN-TV$_1$ &   9.0336     &   8.7566      &  13.3296     &    -    &    7.3971     &   8.5761      &      9.7250 &    10.7706     \\
				\hline TNN-TV$_2$ &   10.1656     &    10.9022     &  12.5970     &   14.1778     &    8.5234     &    9.1056    &      9.6482 &    10.9913    \\
				\hline MF-TV &   50.3907     &    184.2496     &  -     &  -      &    17.7529     &   18.0318      &  -     &   -      \\
				\hline SiLRTC-TT & 71.1641   &  67.1689   &    -   & -       &    15.7843     &    16.8340     &   -    &  - \\ 
				\hline    
			\end{tabular}
			\captionof{table}{The cpu  times required for  TNN, TNN-TV$_1$, TNN-TV$_2$, MF-TV and SiLRTC-TT   for $SR=0.2,\, 0.3$  using  the  "Lena" image.}
			\label{tab 2}
			\end{table}
				\noindent In Table  \ref{tab 2}, we reported the execution times needed to achieve the convergence criterion for each method.   As can be seen from this table, the results obtained by   TV$_1$ and TV$_2$ are faster as compared to the ones obtained by the  other three methods.   \\
	\end{subsection}
	
	\begin{subsection}{Videos}
		In this part we  test the performance of our algorithms on some videos. In our example,  we used the   video of "Suzie" of size $128\times 128\times 150$, and we compared the  obtained results with those of  TNN, SiLRTC-TT and MF-TV. In Figure \ref{fig 4} we gave  the recovered results of one frame  for $SR=0.1$ for the first line and $SR=0.05$ for the second one.
		\begin{figure}[H]
			\centering 
			\hspace*{- 0.5 cm}\begin{tabular}{c}
				\includegraphics[width=1.16\linewidth]{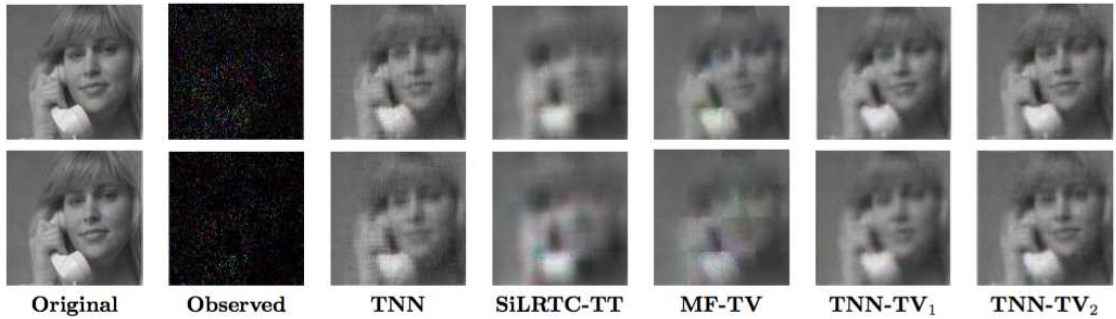}
			\end{tabular} 
			\captionof{figure}{The results of the algorithms TNN, SiLRTC-TT, MF-TV, TNN-TV$_1$ and TNN-TV$_2$ for a frame of the video "Suzie" with $SR=0.1$ for the first line and $SR=0.05$ for the second line} \label{fig 4}
		\end{figure}
		\begin{figure}[H]
			\centering
			\begin{tabular}{cc}
				\includegraphics[width=.5\linewidth]{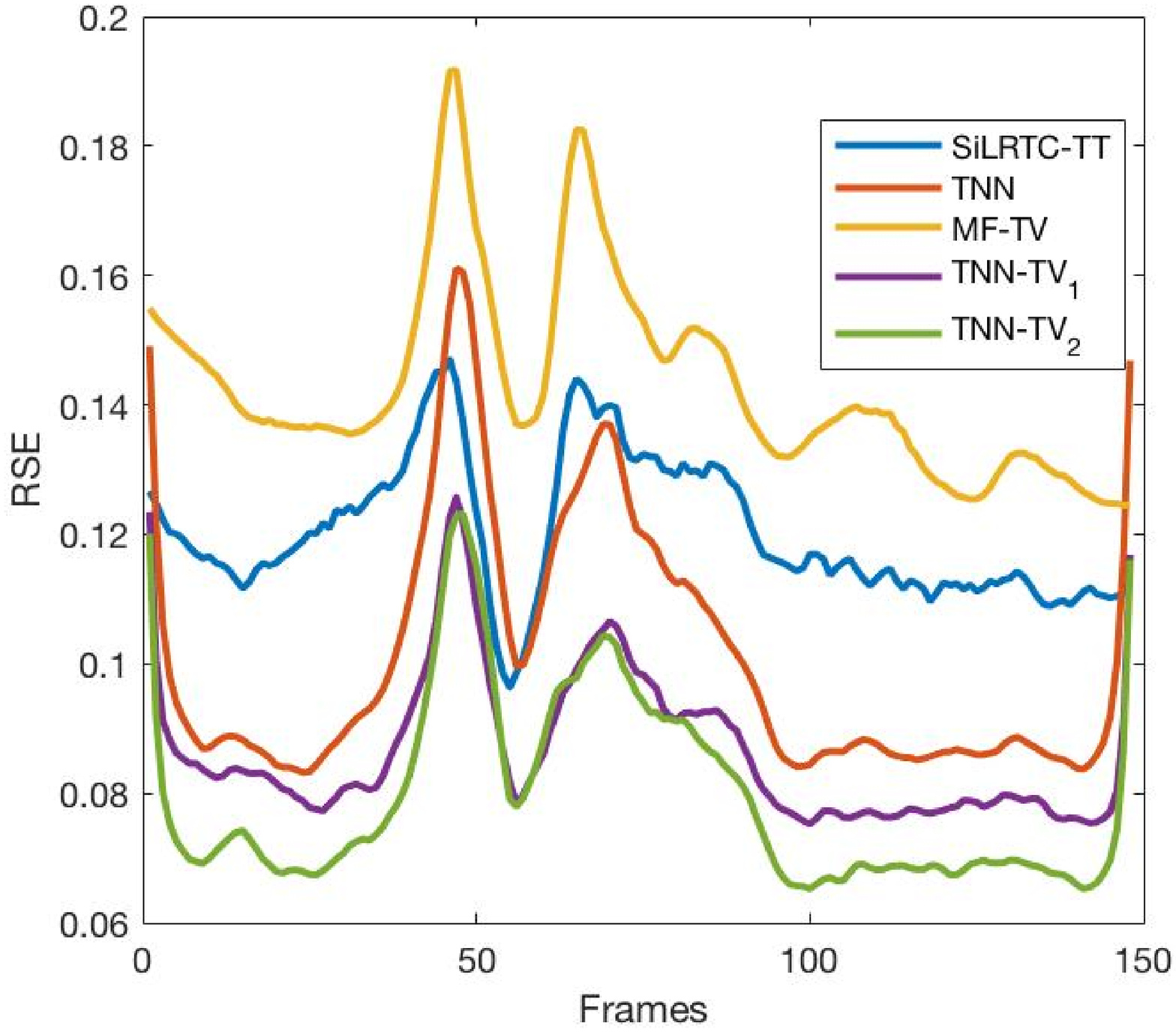}&
				\includegraphics[width=.56\linewidth]{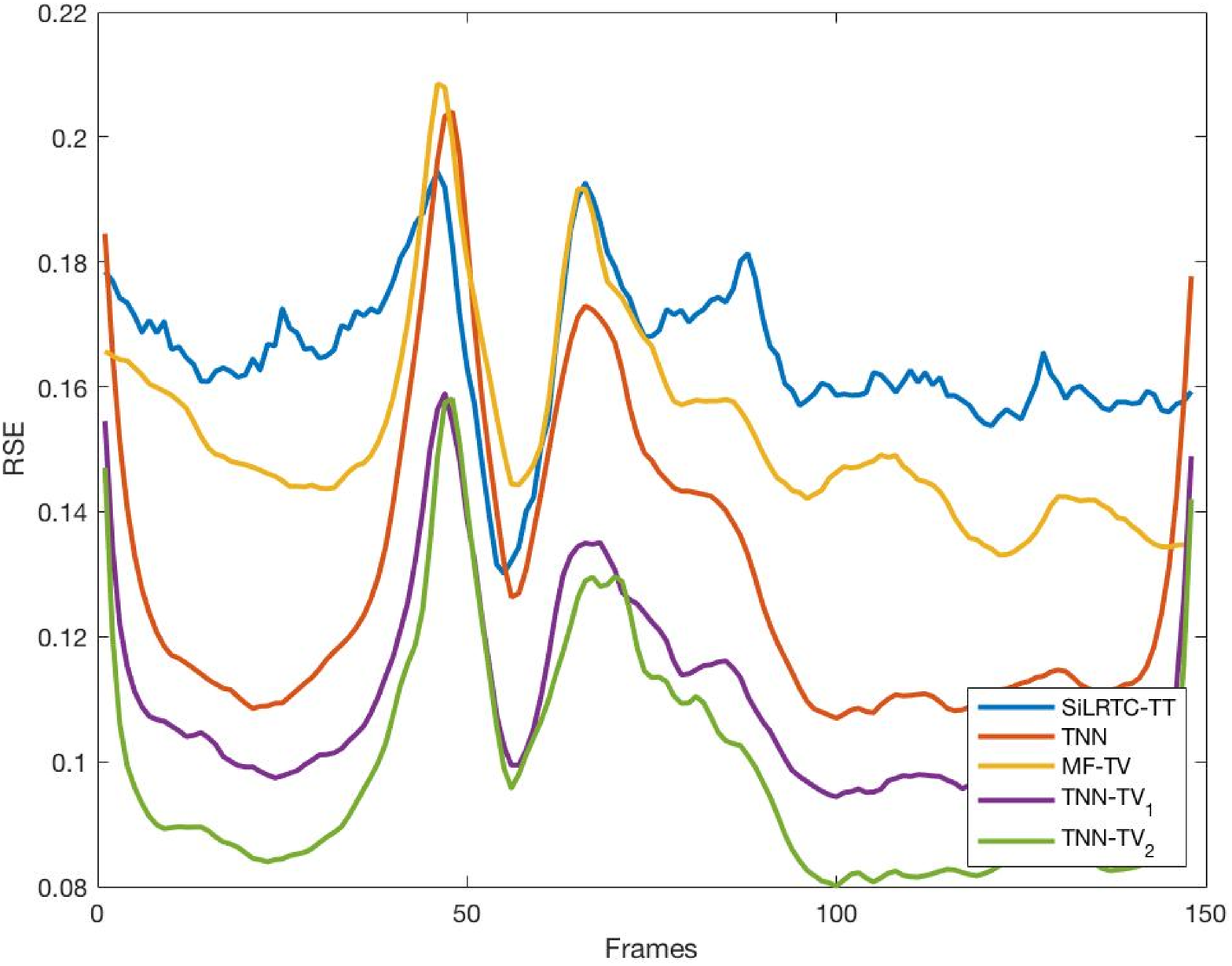}
			\end{tabular}
			\captionof{figure}{The values of RSE by the algorithms SiLRTC-TT, TNN, MF-TV, TNN-TV$_1$ and TNN-TV$_2$ for each frame of the video of "Suzie" for $SR=0.1$ and $SR=0.05$ from the left to the right.} \label{fig 5}
		\end{figure}
		\noindent  In  Figure \ref{fig 5} we plotted the values of RSE for each frame obtained  by  TNN, SiLRTC-TT, MF-TV, TNN-TV$_1$ and TNN-TV$_2$.	As can be seen from this figure,  TNN-TV$_1$ and TNN-TV$_2$ return the best results.  
	\end{subsection}
	
	\begin{subsection}{MRI}
		In this subsection we test our methods  on the MRI data of the front direction.   In this  example, we used  a video of MRI of front direction of size $181 \times  217 \times  150$. In Figure \ref{fig 6} we showed  two recovered frames of this video for $SR=0.1$.  
		\begin{figure}[H]
			\centering
			\begin{tabular}{cccccc}
				\includegraphics[width=1.16\linewidth]{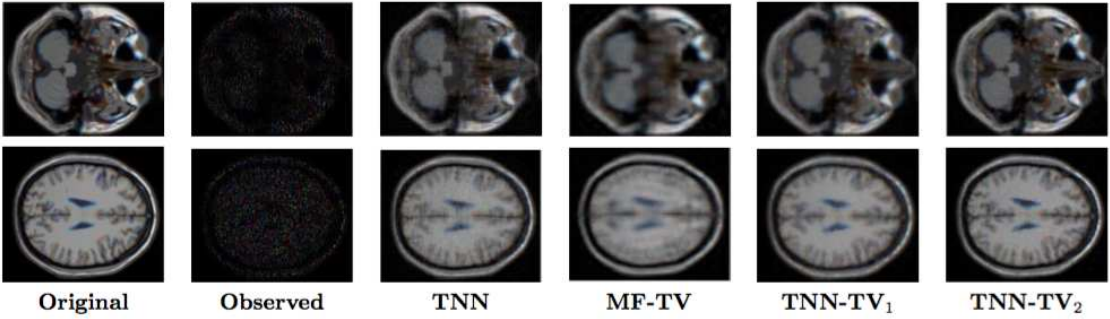}
			\end{tabular}
			\captionof{figure}{The results of the algorithms TNN, MF-TV, TNN-TV$_1$ and TNN-TV$_2$ for the video of the front direction with $SR=0.1$} \label{fig 6}
		\end{figure}
		
		\begin{figure}[H]
			\centering
			\includegraphics[width=.8\linewidth]{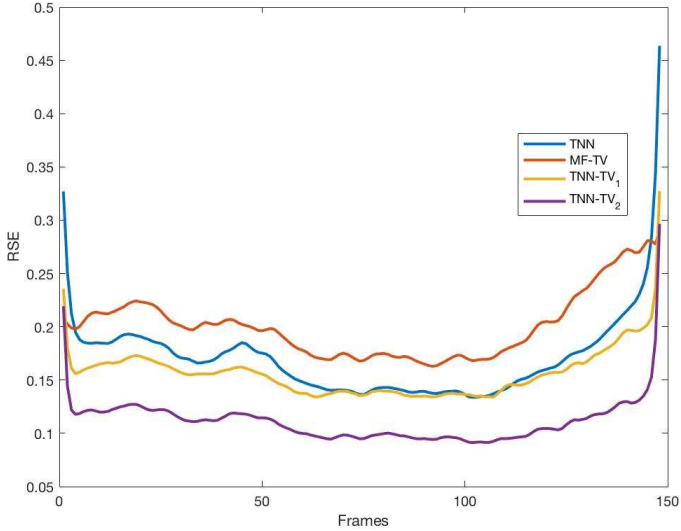}
			\captionof{figure}{The values of RSE  for each frame for the video MRI of the front direction by the algorithms TNN, MF-TV, TNN-TV$_1$ and TNN-TV$_2$.} \label{fig 7}
		\end{figure}
		\noindent In Figure \ref{fig 7} we plotted  the values of RSE for each frame of the video MRI obtained with TNN, MF-TV, TNN-TV$_1$ and TNN-TV$_2$. As one can see from this figure,  TNN-TV$_2$ returns the best result. 
	\end{subsection}
	
	\begin{section}{Conclusion}
		In this paper we proposed two methods for  image completion  by  combining  the tensor nuclear norm and the total variation regularization approaches. We showed how to compute the different tensor sequences obtained from different optimisation problems and gave some convergence theoretical results. The numerical experiments show that our two approaches are efficient and very competitive as compared to other recent completion methods.
	\end{section}

	
\end{document}